\documentclass[rotating,number,citesort,seceqn,dvips]{arxbj}
\usepackage{dcolumn,subenv}
\usepackage{graphicx}

\aid{0}
\volume{17}
\issue{2}
\pubyear{2011}
\firstpage{562}
\lastpage{591}
\doi{10.3150/10-BEJ284}

\makeatletter
\newcommand{\bigtimes}{\mathop{\,\mbox{\parbox[c][9pt][b]{18pt}%
            {\fontsize{18}{18}\selectfont{$\times$}}}\!\!}}
\newtheorem{theorem}{Theorem}[section]
\newtheorem{lemma}[theorem]{Lemma}

\newcolumntype{d}[1]{D{.}{.}{#1}}
\newremark{remark}{Remark}[section]
\newcommand{\eqref}[1]{(\ref{#1})}
\newcommand{\ud}{\mathrm{d}}
\newcommand{\h}{h}
\newcommand{\RR}{\mathbb{R}}
\newcommand{\I}{\mathcal{S}}
\newcommand{\Q}{\mathcal{Q}}
\newcommand{\PP}{\mathbb{P}}
\newcommand{\veccc}[1]{\mathbf{#1}}

\def\sfrac#1#2{#1/#2}
\def\sklfrac#1#2{(#1/#2)}

\newcommand{\bseriese}[1]{#1}
\newcommand{\beditore}[1]{#1}
\newcommand{\bbooktitlee}[1]{\textit{#1}}

\newcommand{\bauthore}[1]{#1}
\newcommand{\bsnme}[1]{#1}
\newcommand{\bfnme}[1]{}
\newcommand{\binitse}[1]{#1}
\newcommand{\byeare}[1]{#1}
\newcommand{\bjournale}[1]{\textit{#1}}
\newcommand{\btitlee}[1]{#1}
\newcommand{\bvolumee}[1]{\textbf{#1}}
\newcommand{\bpagese}[1]{#1}
\newcommand{\baddresse}[1]{#1}
\newcommand{\bpublishere}[1]{#1}
\newcommand{\bnotee}[1]{#1}
\newcommand{\AND}{and }

\makeatother

\begin{document}
\begin{frontmatter}

\title{The AEP algorithm for the fast computation of the distribution
of the sum of dependent random variables}
\runtitle{The AEP algorithm}

\runauthor{P. Arbenz, P. Embrechts and G. Puccetti}

\begin{aug}
\author[a]{\fnms{Philipp} \snm{Arbenz}\thanksref{a,e1}\ead[label=e1,mark]{philipp.arbenz@math.ethz.ch}},
\author[a]{\fnms{Paul} \snm{Embrechts}\thanksref{a,e2}\ead[label=e2,mark]{embrechts@math.ethz.ch}}
\and
\author[b]{\fnms{Giovanni} \snm{Puccetti}\corref{}\thanksref{b}\ead[label=e3]{giovanni.puccetti@unifi.it}}

\runauthor{P. Arbenz, P. Embrechts and G. Puccetti}

\address[a]{Department of Mathematics,
ETH Zurich,
R\"{a}mistrasse 101,
8092 Zurich,
Switzerland.\\
\printead{e1,e2}}

\address[b]{Dipartimento di Matematica per le Decisioni,
via Lombroso 6/17,
50134 Firenze,
Italy.\\
\printead{e3}}
\end{aug}

\received{\smonth{5} \syear{2009}}
\revised{\smonth{1} \syear{2010}}

%
\begin{abstract}
We propose a new algorithm to compute numerically the distribution
function of the sum of $d$ dependent, non-negative random variables
with given joint distribution.
\end{abstract}

%

\end{frontmatter}

\section{Motivations and preliminaries}\label{se:intro}
In probability theory, the exact calculation of the distribution
function of the sum of $d$ dependent random variables $X_1,\ldots,X_d$
is a rather onerous task.
Even assuming the knowledge of the joint distribution $H$ of the vector
$(X_1,\ldots,X_d)$, one often has to rely on tools like Monte Carlo and
quasi-Monte Carlo methods. All of these techniques warrant considerable
expertise and, more importantly, need to be tailored to the specific
problem under consideration.
In this paper, we introduce a numerical procedure, called the AEP
algorithm (from the names of the authors), which accurately calculates
%
\begin{equation}\label{eq:mainissue0}
\PP[X_1+\cdots+X_d \leq s]
\end{equation}
at a fixed real threshold $s$ and only uses the joint distribution $H$
without the need for any specific adaptation.

Problems like the computation of~\eqref{eq:mainissue0} arise especially
in insurance or finance when one has to calculate an overall capital
charge in order to offset the risk position $S_d=X_1+\cdots+X_d$
deriving from a portfolio of $d$ random losses with known joint
distribution $H$. The minimum capital requirement associated to $S_d$
is typically calculated as the value-at-risk (i.e., quantile) for the
distribution of $S_d$ at some high level of probability.
Therefore, the calculation of a VaR-based capital requirement is
equivalent to the computation of the distribution of $S_d$ (see~\eqref
{eq:mainissue0}). For an internationally active bank, this latter task
is required, for example, under the terms of the New Basel Capital
Accord (Basel II); see~\citep{BASEL2d}.

An area of application in quantitative risk management where our
algorithm may be particularly useful is stress-testing. In this
context, one often has information on the marginal distributions of the
underlying risks, but wants to stress-test the interdependence between
these risks; a concept that enters here is that of the copula.
Especially in the context of the current (credit) crisis, flexibility
of the copula used when linking marginal distributions to a joint
distribution has no doubt gained importance; see, for instance,~\citep{pE09b}.

Although the examples treated in this paper are mainly illustrative,
the dimension $d$ ($\leq$5), the marginal assumptions and the
dependence structure (Clayton and Gumbel copula) used are typical for
risk management applications in insurance and finance.
For more information on this type of question, see, for instance,~\citep{SCOR08,ADO07,BDI08}.

In the following, we will denote (row) vectors in boldface, for
example, $\veccc{1}=(1,\ldots,1) \in\mathbb{R}^d$, $d>1$. $\veccc{e}_k$
represents the $k$th vector of the canonical basis of $\mathbb{R}^d$
and $D=\{1,\ldots,d\}$.
Given a vector $\veccc{b}=(b_1,\ldots,b_d) \in\mathbb{R}^d$ and a real
number $h$, $\Q(\veccc{b},h) \subset\mathbb{R}^d$ denotes the hypercube
defined as
%
\begin{eqnarray}\label{de:cube}
\Q(\veccc{b},h)= \cases{
\displaystyle \bigtimes_{k=1}^d (b_k,b_k+h], &\quad \mbox{if }$h> 0$,\cr
\displaystyle \bigtimes_{k=1}^d (b_k+h,b_k], &\quad \mbox{if }$h< 0$.
}
\end{eqnarray}
For notational purposes, we set $\Q(\veccc{b},0)=\varnothing$. On some
probability space $(\Omega, \mathfrak{A},\PP)$, let the random
variables $X_1,\ldots,X_d$ have joint $d$-variate distribution $H$;
$H$ induces the probability measure~$V_H$ on $\RR^d$ via
\[
V_H \Biggl[\displaystyle \bigtimes_{i=1}^{d}(-\infty,x_i] \Biggr]=H(x_1,\ldots,x_d).
\]
We denote by $\veccc{i}_0,\ldots, \veccc{i}_{N}$ all the $2^d$ vectors in
$\{0,1\}^d$, that is, $\veccc{i}_0=(0,\ldots,0)$, $\veccc{i}_k=\veccc{e}_k,
k\in D$, and so on, $\veccc{i}_{N}=\veccc{1}=(1,\ldots,1)$, where
$N=2^{d}-1$. By $\# \veccc{i}=\sum_{k=1}^d i_k$, we denote the number of
$1$'s in the vector $\veccc{i}$, for example, $\# \veccc{i}_0=0, \# \veccc{i}_N=d$.
The $V_H$-measure of a hypercube $ \Q(\veccc{b},h)$, $h>0$, can also be
expressed as
%
\begin{equation}\label{eq:rect}
V_H [ \Q(\veccc{b},h) ]= \mathbb{P} \bigl[X_k \in(b_k,b_k+h], k
\in D  \bigr]=\sum_{j=0}^{N} (-1)^{d- \#\veccc{i}_j} H  ( \veccc{b}
+h\veccc{i}_j  ).
\end{equation}
The case $h<0$ is analogous.
If necessary, \eqref{eq:rect} can also be expressed in terms of the
survival function $\overline{H}=1- H$.
Moreover, $\I(\veccc{b},h) \subset\mathbb{R}^d$ denotes the
$d$-dimensional simplex defined as
%
\begin{eqnarray}\label{de:simplex}
\I(\veccc{b},h) = \cases{
 \Biggl\{\veccc{x} \in\mathbb{R}^d \dvtx  x_k-b_k >0, k\in D \mbox{ and } \displaystyle \sum_{k=1}^d (x_k-b_k) \leq h \Biggr\}, &\quad\mbox{if }$h>0$,\cr
\Biggl\{\veccc{x}\in\mathbb{R}^d \dvtx  x_k-b_k \leq0, k \in D \mbox{ and }\displaystyle  \sum_{k=1}^d (x_k-b_k) > h \Biggr\},  &\quad\mbox{if }$h< 0$.
}
\end{eqnarray}
Again, $\I(\veccc{b},0)=\varnothing$.
Finally, we denote by $\lambda_d$ the Lebesgue measure on $\RR^d$. For
instance, the Lebesgue measure of the simplex $ \I(\veccc{b},h)$ is given by
%
\begin{equation}\label{eq:volumesim}
\lambda_d [\I(\veccc{b},h) ]=\frac{ | h |^d}{d!}.
\end{equation}

\section{Description of the AEP algorithm for $d=2$}\label{se:algo2}
Throughout the paper, we assume the random variables $X_1,\ldots,X_d$
to be non-negative, that is, $\PP[X_k \leq0]=0, k \in D$. The
extension to random variables bounded from below is straightforward and
will be illustrated below.
We assume that we know the joint distribution $H$ of the vector
$(X_1,\ldots,X_d)$ and define $S_d=X_1+\cdots+X_d$. Our aim is then to
numerically calculate
\begin{eqnarray*}
\PP[ S_d \leq s]= V_H [\I(\veccc{0},s) ]
\end{eqnarray*}
for a fixed positive threshold $s$.

\begin{figure}[b]

\includegraphics{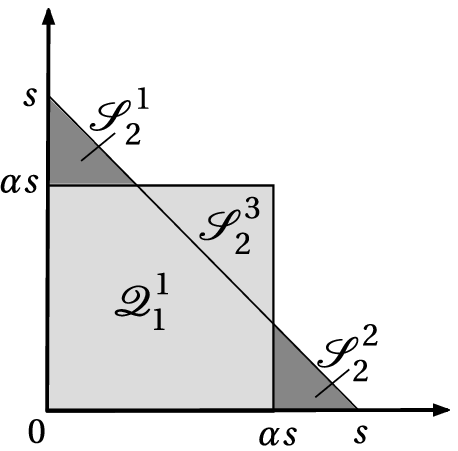}

\caption{Decomposition~\protect\eqref{eq:setdec2beg} of the two-dimensional
simplex $\I(\veccc{0},s)$.}\label{fi:dec2dbeg}
\end{figure}

Due to~\eqref{eq:rect}, it is very easy to compute the $V_H$-measure of
hypercubes in $\mathbb{R}^d$. The idea behind the AEP algorithm is then
to approximate the simplex $\I(\veccc{0},s)$ by hypercubes.
Before proceeding to the general case, we first illustrate our method
for dimension $d=2$.

As illustrated in Figure~\ref{fi:dec2dbeg}, the $V_H$-measure of the
simplex $\I_1^1=\I(\veccc{0},s)$ can be proxied by the $V_H$-measure of
the hypercube $\Q_1^1=\Q(\veccc{0},\alpha s)$ with $\alpha\in[1/2,1)$.
The error committed by using this approximation can be expressed in
terms of the measure of the three simplexes
\begin{eqnarray*}
\I_2^1&=&\I\bigl((0,\alpha s),(1-\alpha)s\bigr),\qquad
\I_2^2=\I\bigl((\alpha s,0),(1-\alpha)s\bigr)\quad \mbox{and}\quad
\\
\I_2^3&=&\I\bigl((\alpha s,\alpha s),(1-2\alpha)s\bigr).
\end{eqnarray*}
Formally, we have
%
\begin{equation}\label{eq:setdec2beg}
\I(\veccc{0},s)= (\Q_1^1 \cup\I_2^1 \cup\I_2^2 )
\setminus\I_2^3\qquad \mbox{for all } \alpha\in[1/2,1).
\end{equation}
Since $\alpha\in[1/2,1)$, the sets $\I_2^1,\I_2^2$ and $\Q_1^1$ are
pairwise disjoint. Also, note that $\I_2^3 \subset\Q_1^1$. The
$V_H$-measure of $\I(\veccc{0},s)$ can thus be written as
\begin{eqnarray*}
V_H [ \I(\veccc{0},s) ]=
V_H [\Q_1^1 ] + V_H [\I_2^1 ]+V_H [\I_2^2 ]-V_H [\I_2^3 ].
\end{eqnarray*}
With the notation $s_2^1=s_2^2=1$ and $s_2^3=-1$, we translate the
equation above into
%
\begin{equation}\label{eq:proof3}
V_H [ \I(\veccc{0},s) ]=
V_H [\Q_1^1 ]+ \sum_{k=1}^3 s_2^k V_H [\I_2^k ].
\end{equation}
Using~\eqref{eq:rect}, a first approximation of $V_H [\I(\veccc{0},s) ]$
is given by the value
\[
P_1(s)=V_H [\Q_1^1 ]=H(\alpha s,\alpha s)-H(0,\alpha s)-H(\alpha s,0)+H(0,0).
\]
Using~\eqref{eq:proof3}, the error committed by considering $P_1(s)$
instead of $V_H [\I(\veccc{0},s) ]$ can be expressed in terms of the
$V_H$-measure of the three simplexes $\I_2^k$ defined above, that is,
%
\begin{equation}\label{eq:n=1}
V_H [\I(\veccc{0},s) ]-P_1(s)= \sum_{k=1}^3 s_2^k V_H [\I_2^k ].
\end{equation}
At this point, we can apply to each of the $\I_2^k$'s a decomposition
analogous to the one given in~\eqref{eq:proof3} for $\I_1^1=\I(\veccc
{0},s)$, in order to obtain a better approximation of their measures
and hence of the measure of $\I_1^1$.
The only difference between the first and the following step is that we
have to keep track of whether the measure of a simplex has to be added
to or subtracted from the next approximation, $P_{2}(s)$, of $V_H [ \I
(\veccc{0},s) ]$. The value $s_2^k$, associated to each simplex $\I_2^k$,
indicates whether the corresponding measure is to be added ($s_2^k=1$)
or subtracted ($s_2^k=-1$).
The next approximation, $P_2(s)$, will be defined such that the
difference $V_H [\I(\veccc{0},s) ]-P_2(s)$ is the sum of the
$V_H$-measures of a total of nine simplexes produced by the
decompositions of the three $\I_2^k$'s. The nine simplexes are then
passed as input to the third iteration and so on.

Before formally defining the algorithm in arbitrary dimension $d$, it
is important to make the following points.
\begin{itemize}
\item We will prove that the set decomposition~\eqref{eq:setdec2beg}
holds analogously in arbitrary dimension $d$ for every choice of
$\alpha\in[1/d,1)$. Unfortunately, the simplexes $\I_{n+1}^k$
generated at the $n$th iteration of the algorithm
are, in general, not disjoint for $d > 2$. This will imply a more
complicated formula for the general $V_H$-measure decomposition.
\item Equation~\eqref{eq:proof3} depends on the choice of $\alpha$. In
Section~\ref{se:gammachoice}, we will study an optimal choice for
$\alpha$.
\end{itemize}

\section{Description of the AEP algorithm for arbitrary $d$}\label{se:algo}
Recall that in Section~\ref{se:intro}, we denoted by $\veccc{i}_0,\ldots,
\veccc{i}_{N}$ all of the $2^d$ vectors in $\{0,1\}^d$, where $N=2^d -1$.
Also, let $\alpha\in[1/d,1)$. At the beginning of the $n$th iteration
($n \in\mathbb{N}$), the algorithm receives as input $N^{n-1}$
simplexes which we denote by $\I_n^k=\I(\veccc{b}_n^k,\h_n^k),$ for $k=1,
\ldots, N^{n-1}$. To each simplex, we associate the value $s_n^k \in\{
-1,1 \}$, which indicates whether the measure of the simplex has to be
added ($s_n^k=1$) or subtracted ($s_n^k=-1$) in order to compute an
approximation of $V_H [\I(\veccc{0},s) ]$.

Each simplex $\I_n^k$ is then decomposed via one hypercube $\Q_n^k=\Q
(\veccc{b}_{n}^k,\alpha\h_{n}^k)$ and $N$ simplexes $\I_{n+1}^k=\I(\veccc
{b}_{n+1}^k, \h_{n+1}^k)$. In Appendix~\ref{se:appendixa}, we prove the
rather technical result that
the $V_H$-measure of
each simplex $\I_n^k$ can be calculated as
%
\begin{eqnarray}\label{eq:decd}
V_H [\I_n^k ]=V_H [\Q_n^k ]+\sum_{j=1}^{N} m^{j} V_H [\I_{n+1}^{Nk-N+j} ],
\end{eqnarray}
where the sequences $\veccc{b}_n^k,\h_n^k$ and $m^j$ are defined by their
initial values $\veccc{b}^1_1=\veccc{0},\h_1^1= s$ and
\begin{eqnarray}\label{eq:defofbd}
\veccc{b}_{n+1}^{Nk-N+j}&=&\veccc{b}_n^k +\alpha\h_n^k \veccc{i}_j,\qquad
\h_{n+1}^{Nk-N+j}=(1- \# \veccc{i}_j \alpha)\h_{n}^{k}, \nonumber
\\[-8pt]\\[-8pt]
m^{j}&=&\cases{
(-1)^{1+\#\veccc{i}_j},   &\quad\mbox{if }$\# \veccc{i}_j < 1/ \alpha$,\cr
0,                        &\quad\mbox{if }$\# \veccc{i}_j = 1/ \alpha$,\cr
(-1)^{d+1-\#\veccc{i}_j}, &\quad\mbox{if }$\# \veccc{i}_j > 1/ \alpha$
}\nonumber
\end{eqnarray}
for all $j=1,\ldots,N$ and $k=1,\ldots,N^{n-1}$.
At this point, we note that by changing the value $\veccc{b}_1^1$, one
can apply the algorithm to the case in which the random vector
$(X_1,\ldots,X_d)$ also assumes negative values, but is still bounded
from below by $\veccc{b}_1^1$.

We define the sequence $P_n(s)$ as the sum of the $V_H$-measures of the
$\Q_n^k$, multiplied by the corresponding $s_n^k$,
%
\begin{equation}\label{eq:defpnd}
P_{n}(s)=P_{n-1}(s)+ \sum_{k=1}^{N^{n-1}} s_n^k V_H [\Q_n^k ]=\sum
_{i=1}^{n}\sum_{k=1}^{N^{i-1}} s_i^k V_H [\Q_i^k ],
\end{equation}
where $P_0(s)=0$ and the $s_n^k$ are defined by $s_1^1=1$ and
%
\begin{equation}\label{eq:defofs}
s_{n+1}^{Nk-N+j}=s_n^km^{j}\qquad \mbox{for all } j=1,\ldots,N \mbox{ and }
k=1,\ldots,N^{n-1}.
\end{equation}

We will show that, under weak assumptions on $H$, the sequence $P_n(s)$
converges to $V_H [\I(\veccc{0},s) ]$.
Moreover, from~\eqref{eq:rect}, $P_n(s)$ can be calculated in a
straightforward way.
The $(N^{n-1})\times N=N^n$ simplexes $\I_{n+1}^k$ generated by~\eqref
{eq:decd} are then passed to the $(n+1)$th iteration in order to
approximate their $V_H$-measures with the measures of the hypercubes $\Q
_{n+1}^k$.

As a first step to show that $P_n(s)$ tends to $V_H [\I(\veccc{0},s) ]$,
we calculate the error by using $P_n(s)$ instead of $V_H [\I(\veccc{0},s) ]$.
\begin{theorem}\label{th:n=nd}
With the notation introduced above, we have that
%
\begin{equation}\label{eq:n=nd}
V_H [\I(\veccc{0},s) ]-P_n(s)=\sum_{k=1}^{N^{n}} s_{n+1}^k V_H [\I_{n+1}^k ]\qquad \mbox{for all } n \in\mathbb{N}.
\end{equation}
\end{theorem}

\begin{pf}
We prove the theorem by induction on $n$. Note that for $n=1$, \eqref
{eq:n=nd} corresponds to~\eqref{eq:decd}.
Now, assume by induction that
\begin{eqnarray*}
V_H [\I(\veccc{0},s) ]=P_{n-1}(s)+\sum_{k=1}^{N^{n-1}} s_{n}^k V_H [\I
_{n}^k ],
\end{eqnarray*}
which, recalling~\eqref{eq:decd},~\eqref{eq:defpnd} and~\eqref
{eq:defofs}, yields
\begin{eqnarray*}
V_H [\I(\veccc{0},s) ]&=&P_{n-1}(s) +\sum_{k=1}^{N^{n-1}} s_{n}^k V_H [\Q
_{n}^k ]+\sum_{k=1}^{N^{n-1}} s_{n}^k \Biggl(\sum_{j=1}^{N} m^{j} V_H [\I
_{n+1}^{Nk-N+j} ] \Biggr) \\
&=& P_n(s) + \sum_{k=1}^{N^{n-1}} \sum_{j=1}^{N} s_{n}^k m^{j} V_H [\I
_{n+1}^{Nk-N+j} ]\\
&=&P_n(s) + \sum_{k=1}^{N^{n-1}} \sum_{j=1}^{N} s_{n+1}^{Nk-N+j} V_H [\I
_{n+1}^{Nk-N+j} ]=P_n(s) + \sum_{k=1}^{N^n} s_{n+1}^{k} V_H [\I
_{n+1}^{k} ].
\end{eqnarray*}
\upqed
\end{pf}

We are now ready to give a sufficient condition for the convergence of
the sequence $P_n(s)$ to $V_H [\I(\veccc{0},s) ]$.
The idea of the proof is that if the total Lebesgue measure of the new
$N$ simplexes $\I_{n+1}^{Nk-N+j},j=1,\ldots,N,$ generated by the
simplex $\I_n^k$, is smaller than the Lebesgue measure of $\I_n^k$
itself, then, by assuming continuity of $H$, the error~\eqref{eq:n=nd}
will go to zero.
Let us define $e_n=\sum_{k=1}^{N^n} \lambda_d [\I_{n+1}^k ]$ to be the
sum of the Lebesgue measure of the simplexes passed to iteration $n+1$.
We define the \textit{volume factor} $f(\alpha)$
to be the ratio between the sum of the Lebesgue measure of the
simplexes in two subsequent iterations, that is,
$f(\alpha)=e_n/e_{n-1}$. Recalling the formula~\eqref{eq:volumesim} for
the $\lambda_d$-measure of a simplex, we have that
\begin{eqnarray*}
\sum_{j=1}^N \lambda_d [\I_{n+1}^{Nk-N+j} ]=\sum_{j=1}^N \frac{ | ( 1-\#
\veccc{i}_j \alpha ) \h_n^k |^d}{d!}=\sum_{j=1}^d \pmatrix{d\cr j} \frac{
| 1-j \alpha |^d |\h_n^k |^d}{d!}.
\end{eqnarray*}
Observing that
the $N$ simplexes $\I_{n+1}^{Nk-N+j}$, $j=1,\ldots,N$,$d\choose j$ are generated by
the simplex $\I_n^k$,
we use the above equation to conclude that
\begin{eqnarray*}
f(\alpha)&=&\frac{e_n}{e_{n-1}}=\frac{\sum_{k=1}^{N^n} \lambda_d [\I
_{n+1}^k ]}{\sum_{k=1}^{N^{n-1}} \lambda_d [\I_{n}^k ]}\\
&=&\frac{\sum
_{k=1}^{N^{n-1}}\sum_{j=1}^N \lambda_d [\I_{n+1}^{Nk-N+j} ]}{\sum
_{k=1}^{N^{n-1}} \lambda_d [\I_{n}^k ]}=\frac{\sum_{k=1}^{N^{n-1}} \sum
_{j=1}^d \left({d\atop j}\right) \sklfrac{ | 1-j \alpha |^d |\h_n^k |^d}{d!} }{\sum
_{k=1}^{N^{n-1}} \lambda_d [\I_{n}^k ]}
\\
&=&\frac{\sklfrac{1}{d!}\sum_{k=1}^{N^{n-1}} |\h_n^k |^d \sum_{j=1}^d
\left({d\atop j}\right) | 1-j \alpha |^d}{\sklfrac{1}{d!}\sum_{k=1}^{N^{n-1}} |\h
_n^k |^d}
=\sum_{j=1}^d \pmatrix{d\cr j} | 1-j \alpha |^d.
\end{eqnarray*}

A sufficient condition for the convergence of the AEP algorithm can
then be expressed in terms of the volume factor $f(\alpha)$. We first
assume $H$ to be absolutely continuous with a bounded density.
\begin{theorem}\label{th:convd}
Assume that $V_H$ has a bounded density $v_H$. If the volume factor
satisfies \mbox{$f(\alpha)<1,$}
then
%
\begin{equation}\label{eq:maintheorem}
\lim_{n \to\infty} P_n(s) = V_H [\I(\veccc{0},s) ].
\end{equation}
\end{theorem}

\begin{pf}
Since $V_H$ has a density $v_H$ bounded by a constant $c>0$,
using~\eqref{eq:n=nd}, we have that
\begin{eqnarray*}
| V_H [\I(\veccc{0},s) ]-P_n(s) |
&=& \Bigg|\sum_{k=1}^{N^{n}} s_{n+1}^k V_H [\I_{n+1}^k ] \Bigg|
= \Bigg|\sum_{k=1}^{N^{n}} \int_{\I_{n+1}^k} s_{n+1}^k \,\ud H \Bigg|
\\& \leq& \sum_{k=1}^{N^{n}}
\bigg| \int_{\I_{n+1}^k} s_{n+1}^k c \,\ud\lambda_d \bigg| \leq c\sum_{k=1}^{N^{n}} \int_{\I_{n+1}^k} | s_{n+1}^k | \,\ud
\lambda_d
\\&=& c \sum_{k=1}^{N^n} \int_{\I_{n+1}^k} \ud\lambda_d =
c \sum_{k=1}^{N^n} \lambda_d [\I_{n+1}^k ]=c e_n.
\end{eqnarray*}
We conclude by noting that since $e_n >0$ and $e_{n}/e_{n-1}=f(\alpha
)<1$ by assumption, $e_n$ goes to zero exponentially in $n$.
\end{pf}

In order for~\eqref{eq:maintheorem} to hold, it is sufficient that
$v_H$ is bounded on $\bigcup_{k=1}^{N^n} \I_{n+1}^k$ for $n$ large enough.
Define the curve $\Gamma_s$ as
%
\begin{equation}\label{eq:gammasdef}
\Gamma_s= \Biggl\{(x_1,\ldots,x_d) \in\mathbb{R}^d\dvtx  \sum_{k=1}^{d} x_k=s \Biggr\}.
\end{equation}
The following theorem states that the $L^1$-distance from the curve
$\Gamma_s$ to each point in $\bigcup_{k=1}^{N^n} \I_{n+1}^k$ is bounded
by a factor $\gamma^n s$, where $\gamma=\max \{1-\alpha, | 1-d \alpha
| \}$. When $\alpha\in(0,2/d)$, we have $\gamma<1$ and that this
distance goes to zero as $n \to\infty$. For Theorem~\ref{th:convd} to
hold when $\alpha\in(0,2/d)$, it is then sufficient to require that
$H$ has a bounded density only in a neighborhood of $\Gamma_s$.
We will discuss this assumption further in Section~\ref{se:conclusions}.
\begin{theorem}\label{th:interm}
If $\veccc{x} \in\bigcup_{k=1}^{N^n} \I_{n+1}^k$, then its
$L^1$-distance from the curve $\Gamma_s$ is bounded by $\gamma^n s$,
with $\gamma=\max \{1-\alpha, | 1-d \alpha | \}$.
\end{theorem}
%

\begin{pf} We denote by $b_n^{k,r}$ (resp., $i_j^{r}$) for $r \in
D$ the $d$ components of the vectors $\veccc{b}_n^k$ (resp., $\veccc
{i}_j$). We prove by induction on $n$ that
%
\begin{equation}\label{eq:againind}
\sum_{r=1}^{d}b_n^{k,r}+\h_n^k=s\qquad \mbox{for all } k=1,\ldots,N^{n-1}
\mbox{ and } n \geq1.
\end{equation}
For $n=1$, the statement is true since there is only one simplex with
$\veccc{b}_1^{1}=\veccc{0}$ and $\h_1^1=s$.
Now, assume the statement holds for $n>1$. By~\eqref{eq:defofbd}, we
have that, for all $j=1,\ldots,N$ and $k=1,\ldots,N^{n-1}$,
\begin{eqnarray*}
&&\sum_{r=1}^{d} b_{n+1}^{Nk-N+j,r}+\h_{n+1}^{Nk-N+j}
\\
&&\quad=\sum_{r=1}^{d}
(b_n^{k,r}+\alpha\h_n^k i_j^{r} )+(1-\# \veccc{i}_j \alpha)\h_{n}^k=
\sum_{r=1}^{d} b_n^{k,r}+\alpha\h_n^k \sum_{r=1}^{d} i_j^{r}+ \h_{n}^k
-\alpha\h_{n}^k\# \veccc{i}_j
\\
&&\quad=\sum_{r=1}^{d} b_n^{k,r}+\alpha\h_n^k \# \veccc{i}_j+ \h_{n}^k -\alpha\h
_{n}^k \# \veccc{i}_j=\sum_{r=1}^{d} b_n^{k,r}+ \h_{n}^k=s,
\end{eqnarray*}
where the last equality is the induction assumption.
Due to~\eqref{eq:againind}, every simplex $\I_{n+1}^k$ generated by the
AEP algorithm has its diagonal face lying on the curve $\Gamma_s$.
As a consequence, the $L^1$-distance from $\Gamma_s$ of each point in
$\I_{n+1}^k$ is strictly smaller than the distance of the vector $\veccc
{b}_{n+1}^k$, which is $ | \h_{n+1}^k |$.
For a fixed $n$ and $k=1,\ldots N^{n-1}$, we have that $
|h^{Nk-N+j}_{n+1} | \leq\gamma | h^k_{n} |$ for all $j=1,\ldots,N$. Hence,
%
\begin{equation}\label{eq:factorga}
\max_{k=1,\ldots,N^{n}} | \h_{n+1}^k | = \gamma^n \h_1^1=\gamma^n s,
\end{equation}
where, for every $n \geq1$, equality holds since we have $
|h^{Nk-N+j}_{n+1} |= \gamma | h_n^k | $ for $j=1$ or $j=N$.
\end{pf}

\section{Choice of $\alpha$}\label{se:gammachoice}
As already remarked, the AEP algorithm depends on the choice of the
parameter $\alpha$.~It~is important to note that, in general, an
optimal choice of $\alpha$ would depend on the measure $V_H$.
In the proof of Theorem~\ref{th:convd}, we have shown that
\[
| P_n(s)-V_H [ \I(\veccc{0},s) ] | \leq C f(\alpha)^n,
\]
where $C$ is a positive constant.
Since we want to keep our algorithm independent of the choice of the
distribution $H$, we suggest using the $\alpha^*$ which minimizes
$f(\alpha)$, that is,
\[
\alpha^*=\mathop{\mathrm{argmin}}_{\alpha\in [\sfrac{1}{d}, 1 )}f(\alpha)=\frac{2}{d+1}.
\]
For dimensions $d\leq7$, some values of $\alpha^*$ and the
corresponding optimal volume factors $f(\alpha^*)$ are given in
Table~\ref{ta:volfact}.
%
\begin{table}[b]
\tablewidth=8cm
\caption{Values for $\alpha^*$ and $f(\alpha^*)$ for dimensions $d \leq
7$}\label{ta:volfact}
\begin{tabular*}{\tablewidth}{@{\extracolsep{4in minus 4in}}llllll@{}}
\hline
 $d$ & $\alpha^*$ & $f(\alpha^*)$ & $d$ & $\alpha^*$ &$f(\alpha^*)$  \\
\hline
 $2$ & $\frac{2}{3}$ & $\frac{1}{3}$ & $5$ & $\frac{1}{3}$ &$\frac{23}{27}$
 \\[3pt]
 $3$ & $\frac{1}{2}$ & $\frac{1}{2}$ & $6$ & $\frac{2}{7}$ &$>$1  \\[3pt]
 $4$ & $\frac{2}{5}$ & $\frac{83}{125}$ & $7$ & $\frac{1}{4}$ & $>$1  \\
\hline
\end{tabular*}
\end{table}

We will show that using $\alpha^*$ has several desirable consequences.
First, when $\alpha=\alpha^*$ and the dimension $d$ is odd, in the
measure decomposition~\eqref{eq:decd}, a number of
$\bigl({\fontsize{7.6}{7.6}\selectfont\matrix{d\cr(d+1)/2}}\bigr)$ simplexes have the corresponding coefficient $m^j$ equal to
zero and can therefore be neglected, increasing the computational
efficiency of the algorithm.
For example, in the decomposition of a three-dimensional simplex, the
algorithm generates only $4$ new simplexes at every iteration with
$\alpha=\alpha^*$, instead of the $2^d-1=7$ generated with any other
feasible value of $\alpha$. Hence, for $\alpha=\alpha^*$, the number of
new simplexes generated at each step is given by the function
%
\begin{equation}\label{de:fsd}
f_S(d) = \cases{ 2^d-1, &\quad \mbox{if $d$ is even,}\cr
2^d-1-\pmatrix{d\cr (d+1)/2}, &\quad\mbox{if $d$ is odd};
}
\end{equation}
see Section~\ref{se:numerics} for further details on this.

Since we have that (proof of Theorem~\ref{th:interm})
%
\begin{equation}\label{eq:last}
(0,+\infty)^d \cap \Biggl(\bigcup_{k=1}^{N^{n-1}} \I_{n}^k \Biggr) \subset\I\bigl(\veccc
{0},(1+\gamma^n)s\bigr) \big\backslash\I\bigl(\veccc{0},(1-\gamma^n)s\bigr),
\end{equation}
the choice of $\alpha=\alpha^*$ will be convenient.
Note that, when $\alpha=\alpha^* \in(0,2/d)$, we have that $\gamma<1$
and $\gamma^n s$ goes to zero as $n \to\infty$.
In order to guarantee the convergence of the sequence $P_n$, it is then
sufficient to require that the distribution $H$ has a bounded density
only in a neighborhood of $\Gamma_s$.
Moreover, it is straightforward to see that $\alpha^*$ also minimizes
$\gamma$.

As illustrated in Table~\ref{ta:volfact}, Theorem~\ref{th:convd} states
the convergence of the sequence $P_n(s)$ when $d \leq5$.
Various elements affect the speed at which $P_n(s)$ converges. First,
in order to seriously affect the convergence rate of $P_n(s)$, it is,
in general, always possible to put probability mass in a smooth way in
a neighborhood of the curve $\Gamma_s$. For the distributions of\vadjust{\goodbreak}
financial and actuarial interest used in Section~\ref{se:applications},
the algorithm performs very well; slow convergence is typically
restricted to more pathological cases, such as those illustrated in
Section~\ref{se:conclusions}. We also have to consider that, for the
same distribution $H$, it is, in general, required to compute the
distribution of $S_d$ at different thresholds $s$. Problems such as
those described in Section~\ref{se:conclusions} may then occur only at
a few points $s$.

A more relevant issue is the fact that the memory required by the
algorithm to run the $n$th iteration increases exponentially in $n$. At
each iteration of the algorithm, every simplex $\I_n^k$ produces one
hypercube and a number $f_S(d)$ of new simplexes to be passed to the
following iteration; see~\eqref{de:fsd}. The computational effort in
the $(n-1)$th step thus increases as $\mathrm{O}(f_S(d)^n)$.
While the dimensions $d \leq5$ are manageable, as reported in
Section~\ref{se:applications}, the numerical complexity for $d \geq6$
increases considerably and quickly exhausts the memory of a standard computer.

Finally, choosing $\alpha=\alpha^*$ also allows the accuracy of the AEP
algorithm to be increased and, under slightly stronger assumptions on
$H$, will lead to convergence of AEP in higher dimensions; see
Section~\ref{se:numerics}.

%
\begin{figure}

\includegraphics{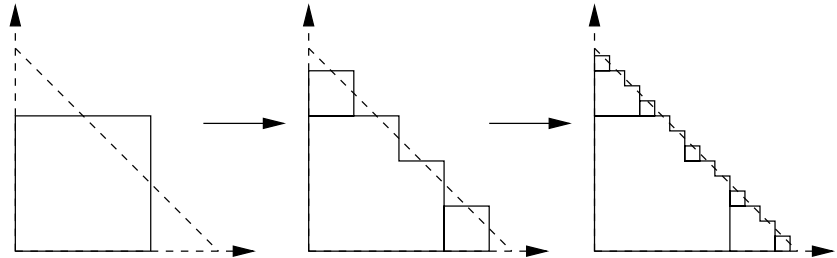}

\caption{The decomposition of a simplex by the AEP algorithm in the case $d=2.$}\label{fi:dec2d}
\vspace*{-5pt}
\end{figure}
\begin{figure}

\includegraphics{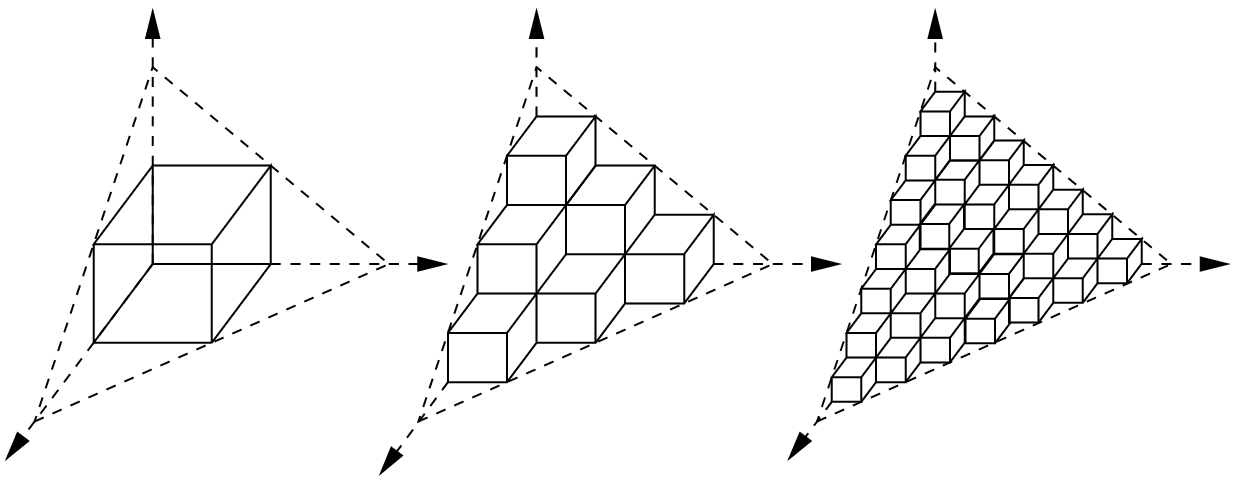}

\caption{The decomposition of a simplex by the AEP algorithm in the case $d=3$.}\label{fi:dec3d}
\end{figure}

We now give some examples of the first step ($n=1$) of the measure
decomposition~\eqref{eq:decd} obtained by choosing $\veccc{b}=\veccc{0}$,
$s=1$, $\alpha=\alpha^*$, for $d=2,3$:
\begin{itemize}
\item in the case $d=2$, with $\alpha= 2/3$, we obtain (see Figure~\ref{fi:dec2d})
\begin{eqnarray*}
V_H [\I((0,0),1) ] &=& V_H \bigl[\Q\bigl((0,0),2/3\bigr) \bigr]+ V_H \bigl[\I\bigl((0,2/3),1/3\bigr) \bigr]
\\
&&{}+ V_H \bigl[\I\bigl((2/3,0),1/3\bigr) \bigr] - V_H \bigl[\I\bigl((2/3,2/3),-1/3\bigr) \bigr];
\end{eqnarray*}
\item in the case $d=3$, with $\alpha= 1/2$, we obtain (see
Figure~\ref{fi:dec3d})
\begin{eqnarray*}
V_H \bigl[\I\bigl((0,0,0),1\bigr) \bigr]&=& V_H \bigl[\Q\bigl((0,0,0),1/2\bigr) \bigr] + V_H \bigl[\I\bigl((1/2,0,0),1/2\bigr) \bigr]
\\
&&{}+ V_H \bigl[\I\bigl((0,1/2,0),1/2\bigr) \bigr] + V_H \bigl[\I\bigl((0,0,1/2),1/2\bigr) \bigr]
\\
&&{}- V_H \bigl[\I\bigl((1/2,1/2,1/2),-1/2\bigr) \bigr].
\end{eqnarray*}
%
\end{itemize}

\section{An improvement of the numerical accuracy of the algorithm via
extrapolation} \label{se:numerics}

In this section, we introduce a method to increase the accuracy of the
AEP algorithm. This method is based on the choice $\alpha=\alpha^*$, as
discussed in Section~\ref{se:gammachoice}.
To this end, we will make the stronger assumption that the joint
distribution $H$ has a twice continuously differentiable density $v_H$,
with bounded derivatives. This will allow us to approximate the density
$v_H$ by its linear Taylor expansion, providing a good estimate of the
approximation error of AEP after a number of iterations.

We first need two simple integration results. Denoting by $\I_{d-1}$ a
simplex in dimension $(d-1)$, for all $s>0$, we have
\begin{eqnarray*}
\int_{\I(\veccc{0},s)} x_d \,\ud\veccc{x}
&=&\int_0^s\!\!\int_0^{s-x_d}\ldots  \int
_0^{s-\sum_{k=3}^d x_k} \int_0^{s-\sum_{k=2}^d x_k} x_d \,\ud\veccc{x}\\
&=&\int_0^s x_d \int_0^{s-x_d}\ldots \int_0^{s-\sum_{k=3}^d x_k} \int
_0^{s-\sum_{k=2}^d x_k} \ud\veccc{x}\\
&=&\int_0^s x_d \lambda_{d-1} [\I_{d-1}(\veccc{0},s-x_d) ] \,\ud x_d = \int
_0^s x_d\dfrac{(s-x_d)^{d-1}}{(d-1)!}\,\ud x_d=\frac{s^{d+1}}{(d+1)!}.
\end{eqnarray*}
Analogously, for all $s>0$, we have
\begin{eqnarray*}
\int_{\Q(\veccc{0},\alpha s)} x_d \,\ud\veccc{x}
&=&\int_0^{\alpha s}\!\!\int_0^{\alpha s}\ldots \int_0^{\alpha s}x_d \,\ud\veccc{x}\\
&=&\int_0^{\alpha s}x_d \int_0^{\alpha s}\ldots \int_0^{\alpha s} \ud\veccc{x}
= (\alpha s)^{d-1} \int_0^{\alpha s} x_d\,\ud x_d
=1/2 (\alpha s )^{d+1}.
\end{eqnarray*}

We now compute the $V_H$-measures of a hypercube and a simplex in the
basic case in which the distribution $H$ has a linear density, that is,
$ v_H(\veccc{b}+\veccc{x}) = a+\sum_{k=1}^d c_k x_k$ for $\veccc{x}\in\I(\veccc
{0}, s) \cup\Q(\veccc{0},\alpha s)$.
For all $s>0$, we obtain
\begin{eqnarray}\label{eq:int1}
V_H [\I(\veccc{b},s) ] &=& a \int_{\I(\veccc{0},s)}\ud\veccc{x} +\sum
_{k=1}^d c_k\int_{\I(\veccc{0},s)}x_k\,\ud\veccc{x} \nonumber\\[-8pt]\\[-8pt]
&=& a\frac{s^d}{d!}+\frac{s^{d+1}}{(d+1)!} \Biggl( \sum_{k=1}^d c_k \Biggr)
= \frac{s^d}{d!} \Biggl(a+\frac{s}{d+1}\sum_{k=1}^d c_k \Biggr),\nonumber \\
\label{eq:int2}
V_H (\Q(\veccc{b},\alpha s ) )
&=& a \int_{\Q(\veccc{0},\alpha s )}\ud\veccc{x} +\sum_{k=1}^d c_k \int
_{\Q(\veccc{0},\alpha s )}x_k\,\ud\veccc{x} \nonumber\\[-8pt]\\[-8pt]
&=& a (\alpha s)^d +\frac{1}{2} \Biggl( \sum_{k=1}^d c_k \Biggr)(\alpha s)^{d+1}
=(\alpha s) ^d \Biggl(a+ \frac{1}{2}\alpha s\sum_{k=1}^d c_k \Biggr).\nonumber
\end{eqnarray}

Thus, for a linear density $v_H$, the ratio $V_H [\I(\veccc{b},s) ]/V_H
[\Q(\veccc{b},\alpha s) ]$ can be made independent of the parameters $\veccc
{b},s,a$ and of the $c_k$'s, by choosing $\alpha=\alpha^*=\frac
{2}{d+1}$, for which we have
%
\begin{equation}\label{eq:ratio4}
V_H [\I(\veccc{b},s) ]=\frac{(d+1)^d}{2^d d!}V_H [\Q(\veccc{b},\alpha^* s) ].
\end{equation}
With similar computations, we obtain the same result for $s<0$.
The following theorem shows that~\eqref{eq:ratio4} analogously holds
for any sufficiently smooth density, in the limit as the number $n$ of
iterations of the AEP algorithm goes to infinity.
\begin{theorem}\label{th:ratiolimit}
Assume that $H$ has a twice continuously differentiable density $v_H$
with all partial derivatives of first and second-order bounded by some
constant $D$. We then have that
%
\begin{equation}\label{eq:ratiolimit}
\sup_{n \in\mathbb{N}} \max_{k=1,\ldots,N^{n-1}} \frac{1}{\vert h_n^k
\vert^{d+2}} \bigg| V_H [\I(\veccc{b}_n^k,h_n^k) ]-
\frac{(d+1)^d}{2^d d!}V_H [\Q(\veccc{b}_n^k,\alpha^* h_n^k) ] \bigg| \leq A <
\infty
\end{equation}
for some positive constant $A$ depending only on the dimension $d$ and
the distribution $H$.
\end{theorem}

\begin{pf}
For a given $\mathbf{b}_n^k$, we can use a Taylor expansion to find some
coefficients $a$ and $c_k$, $k=1,\ldots,d$, depending on $\mathbf{b}_n^k$, such that
%
\begin{equation}\label{eq:taylor}
v_H(\veccc{b}_n^k+\veccc{x})=a+\sum_{k=1}^d c_k x_k+\sum_{ | \beta |
=2}R_{\beta}(\veccc{x})\veccc{x}^\beta \qquad\mbox{for all }\veccc{x} \in\mathcal{B}(\veccc{b}_n^k),
\end{equation}
where $\mathcal{B}(\veccc{b}_n^k)$ is a ball in $\mathbb{R}^d$ centered
at $\veccc{b}_n^k$ such that $\mathcal{B}(\veccc{b}_n^k) \supset\I(\veccc
{b}_n^k,h_n^k) \cup\Q(\veccc{b}_n^k,\alpha^* h_n^k) $.
Note that in equation~\eqref{eq:taylor}, we used multi-index notation\vadjust{\goodbreak}
to indicate that the sum in the last equation extends over
multi-indices $\beta\in\mathbb{N}^d$. Using the assumption on the
partial derivatives of $v_H$, the remainder term $R_{\beta}(\veccc{x})$
satisfies the inequality
%
\begin{equation}\label{eq:taylor2}
| R_{\beta}(\veccc{x}) | \leq \sup_{\veccc{x} \in\mathcal{B}(\veccc{b}_n^k)}
\bigg| \frac{1}{\beta!} \frac{\partial^\beta v_H(\veccc{x})}{\partial\veccc{x}^{\beta}} \bigg| \leq D
\end{equation}
for all $\beta$ with $ | \beta |=2$.
Using~\eqref{eq:taylor} and recalling the expressions~\eqref{eq:int1}
and \eqref{eq:int2} for a linear density and a positive $h_n^k$, we obtain
\begin{eqnarray*}
&&\bigg| V_H [\I(\veccc{b}_n^k,h_n^k) ]-
\frac{(d+1)^d}{2^d d!}V_H [\Q(\veccc{b}_n^k,\alpha h_n^k) ] \bigg| \\
&&\quad=\Bigg| \frac { (h_n^k )^d}{d!} \Biggl(a+\frac{h_n^k}{d+1}\sum_{k=1}^d c_k \Biggr) + \int_{\I(\veccc
{0},h_n^k)} \sum_{\vert\beta\vert=2}R_{\beta}(\veccc{x})\veccc{x}^\beta\,\ud\veccc{x} \\
&&\qquad{}-\frac{(d+1)^d}{2^d d!} \Biggl( ( \alpha h_n^k )^d \Biggl(a+ \frac{1}{2}\alpha
h_n^k \sum_{k=1}^d c_k \Biggr) + \int_{\Q(\veccc{0},\alpha h_n^k)} \sum_{\vert
\beta\vert=2}R_{\beta}(\veccc{x})\veccc{x}^\beta\,\ud\veccc{x} \Biggr) \Bigg|.
\end{eqnarray*}
Choosing $\alpha=\alpha^*$, the previous expression simplifies to
\begin{eqnarray*}
&& \bigg| V_H [\I(\veccc{b}_n^k,h_n^k) ]-
\frac{(d+1)^d}{2^d d!}V_H [\Q(\veccc{b}_n^k,\alpha^* h_n^k) ] \bigg| \\
&&\quad= \bigg|\int_{\I(\veccc{0},h_n^k)} \sum_{\vert\beta\vert=2}R_{\beta}(\veccc{x})\veccc{x}^\beta\,\ud\veccc{x}
-\frac{(d+1)^d}{2^d d!} \int_{\Q(\veccc{0},\alpha^* h_n^k)} \sum_{\vert
\beta\vert=2}R_{\beta}(\veccc{x})\veccc{x}^\beta\,\ud\veccc{x} \bigg|\\
&&\quad\leq \bigg| \sum_{\vert\beta\vert=2} \int_{\I(\veccc{0},h_n^k)} R_{\beta}(\veccc{x})\veccc{x}^\beta\,\ud\veccc{x} \bigg|
+ \frac{(d+1)^d}{2^d d!} \bigg| \sum_{\vert\beta\vert=2} \int_{\Q(\veccc
{0},\alpha^* h_n^k)} R_{\beta}(\veccc{x})\veccc{x}^\beta\,\ud\veccc{x} \bigg|\\
&&\quad\leq D \biggl( \bigg| \sum_{\vert\beta\vert=2} \int_{\I(\veccc{0},h_n^k)} \veccc
{x}^\beta\,\ud\veccc{x} \bigg|
+ \frac{(d+1)^d}{2^d d!} \bigg| \sum_{\vert\beta\vert=2} \int_{\Q(\veccc
{0},\alpha^* h_n^k)} \veccc{x}^\beta\,\ud\veccc{x} \bigg| \biggr),
\end{eqnarray*}
where the last inequality follows from~\eqref{eq:taylor2}. Using the
facts that
\begin{eqnarray*}
\sum_{\vert\beta\vert=2}\int_{\I(\veccc{0},s)} \veccc{x}^\beta\,\ud\veccc{x}
&=& \sum_{i=1}^{d} \int_{\I(\veccc{0},s)} x_i^2 \,\ud\veccc{x}
+2 \sum_{1\leq i<j \leq d} \int_{\I(\veccc{0},s)} x_i x_j \,\ud\veccc{x}\\
&=&\frac{2ds^{d+2}}{(d+2)!}+\frac{2d(d-1)s^{d+2}}{(d+2)!}
\\
&=&\frac{2d^2 s^{d+2}}{(d+2)!}
\end{eqnarray*}
and
\begin{eqnarray*}
\sum_{\vert\beta\vert=2}\int_{\Q(\veccc{0},\alpha s)} \veccc{x}^\beta\,\ud\veccc{x}
&=& \sum_{i=1}^{d} \int_{\Q(\veccc{0},\alpha s)} x_i^2 \,\ud\veccc
{x} +2 \sum_{1 \leq i<j \leq d} \int_{\Q(\veccc{0},\alpha s)} x_i x_j\, \ud\veccc{x}\\
&=& \frac{d (\alpha s )^{d+2}}{3}+\frac{2d(d-1) (\alpha s
)^{d+2}}{4}= \frac{d(3d-1) (\alpha s )^{d+2}}{6},
\end{eqnarray*}
we finally obtain
%
\begin{eqnarray}\label{eq:ratiolimit2}
 \bigg| V_H [\I(\veccc{b}_n^k,h_n^k) ]-
\frac{(d+1)^d}{2^d d!}V_H [\Q(\veccc{b}_n^k,\alpha^* h_n^k) ] \bigg| \leq A
|h_n^k |^{d+2},
\end{eqnarray}
where $A$ is a positive constant depending only on the dimension $d$
and the distribution $H$. Note that in~\eqref{eq:ratiolimit2}, we write
$h_n^k$ in absolute value in order to consider the completely analogous
case in which $h_n^k$ is negative. Thus, the theorem easily follows
from~\eqref{eq:ratiolimit2}.
\end{pf}

Equation~\eqref{eq:ratiolimit} gives a local estimator of the mass of
the simplex $\I(\veccc{b}_n^k,h_n^k)$ in terms of the volume of the
corresponding hypercube $\Q(\veccc{b}_n^k,h_n^k)$, which is
straightforward to compute:
%
\begin{equation}\label{eq:estformula}
V_H [\I(\veccc{b}_n^k,h_n^k) ] \approx\frac{(d+1)^d}{2^dd! } V_H \biggl[\Q
\biggl(\veccc{b}_n^k,\frac{2h_n^k}{d+1} \biggr) \biggr].
\end{equation}
In the case where the density $v_H$ is sufficiently smooth,
it is then possible, after a number of iterations of AEP, to estimate
the right-hand side of~\eqref{eq:n=nd} by using the approximation~\eqref
{eq:estformula}. This procedure defines the estimator $P^*_n(s)$ as
%
\begin{equation}\label{eq:estformula2}
P^*_{n}(s)=P_{n-1}(s)+\frac{(d+1)^d}{2^d d!}\sum_{k=1}^{N^{n-1}}
s_{n}^k V_H [\Q_{n}^k ].
\end{equation}
In what follows, the use of $P^*_n(s)$ as an approximation of $V_H[\I
(\veccc{0},s)]$ will be referred to as the \textit{extrapolation technique}.
The following theorem shows that $P^*_{n}(s)$ converges to $V_H [\I(\veccc
{0},s) ]$ faster, and in higher dimensions, than $P_{n}(s)$.
\begin{theorem}\label{th:extrapconv}
Under the assumptions of Theorem~\ref{th:ratiolimit}, we have, for
$d\leq8$, that
\begin{eqnarray*}
\lim_{n \to+\infty} P^*_n(s)=V_H [ \I(\veccc{0},s) ].
\end{eqnarray*}
\end{theorem}

\begin{pf}
Using~\eqref{eq:n=nd} and~\eqref{eq:ratiolimit2}~in the
definition~\eqref{eq:estformula2} of $P_n^*(s)$, we obtain
\begin{eqnarray}\label{eq:final4}
E^*(n)&=& | V_H [\I(\veccc{0},s) ] - P^*_{n}(s) |\nonumber\\
&=& \Bigg| V_H [\I(\veccc{0},s) ] - P_{n-1}(s) -\frac{(d+1)^d}{2^d d!}
\sum_{k=1}^{N^{n-1}} s_n^k V_H [ \Q_n^k ] \Bigg|\nonumber\\[-8pt]\\[-8pt]
&=& \Bigg| \sum_{k=1}^{N^{n-1}} s_n^k V_H [ \I_n^k ] -\frac
{(d+1)^d}{2^d d!} \sum_{k=1}^{N^{n-1}} s_n^k V_H [ \Q_n^k ] \Bigg|\nonumber\\
&\leq& \sum_{k=1}^{N^{n-1}} \bigg| V_H [ \I_n^k ] -\frac
{(d+1)^d}{2^d d!}V_H [ \Q_n^k ] \bigg| \leq A\sum_{k=1}^{N^{n-1}} | h_n^k
|^{d+2}= Ae^*_{n-1},\nonumber
\end{eqnarray}
where, for the positive sequence $e^*_n=\sum_{k=1}^{N^n} | h_{n+1}^k
|^{d+2}$, we have that
\begin{eqnarray*}
\frac{e^*_n}{e^*_{n-1}}&=&\frac{\sum_{k=1}^{N^{n-1}} \sum_{j=1}^{N} |
h_{n+1}^{Nk-N+j} |^{d+2}}{\sum_{k=1}^{N^{n-1}} | h_{n}^k |^{d+2}}
=\frac{\sum_{k=1}^{N^{n-1}} \sum_{j=1}^{d}\left({d\atop j}\right) | 1-j\alpha^* |^{d+2}
| h_{n}^{k} |^{d+2}}{\sum_{k=1}^{N^{n-1}} | h_{n}^k |^{d+2}}
\\
&=&\frac
{\sum_{k=1}^{N^{n-1}} | h_{n}^{k} |^{d+2} \sum_{j=1}^{d}\left({d\atop j}\right)|
1-j\alpha^* |^{d+2}}{\sum_{k=1}^{N^{n-1}} | h_{n}^k |^{d+2}}=\sum
_{j=1}^{d}\pmatrix{d\cr j} | 1-j\alpha^* |^{d+2}.
\end{eqnarray*}
The theorem follows by noting that the factor $f_*(d),$ defined as
%
\begin{equation}\label{eq:extrapfactor}
f_*(d)=\sum_{j=1}^{d}\pmatrix{d\cr j} | 1-j\alpha^* |^{d+2}
\end{equation}
is less than 1 for $d \leq8$; see Table~\ref{ta:simplexes}. In these
dimensions, $e_n^*$, and hence $E^*(n)$, converge to zero.
\end{pf}

%
\begin{table}[b]
\caption{Extrapolation error ratio $ f_*(d)$ as defined in~\protect\eqref{eq:extrapfactor},
number $ f_S(d)$ of new simplexes produced at each
iteration and convergence rates of the AEP extrapolation error as a
function of the number of evaluations performed by the algorithm; for
$d=9$, convergence of AEP is not assured (na)}\label{ta:simplexes}
\begin{tabular*}{\textwidth}{@{\extracolsep{4in minus 4in}}ld{2.4}d{2.4} d{2.4}d{2.4}d{2.4} d{2.4}d{3.4}d{3.0}@{\hspace*{-2pt}}}
\hline
 $d$ & 2 & 3 & 4 & 5 & 6 & 7 & 8 & 9  \\ \hline
 $f_*(d)$ & 0.0370 & 0.1250 & 0.2339 &0.3580 &0.4982 &0.6556 &0.8314 &{>}1
 \\[3pt]
 $f_S(d)$ & 3 & 4 & 15 &21 &63 &92 &255 &385  \\[3pt]
 $\frac{\ln f_*(d)}{\ln f_S(d)}$ & -3 & -1.5 & -0.54 &-0.34 &-0.17 &-0.09 &-0.033 &\multicolumn{1}{c@{}}{na}  \\
\hline
\end{tabular*}
\end{table}

We should point out that, due to Theorem~\ref{th:interm}, Theorem~\ref{th:extrapconv} also remains valid in the case where $H$ satisfies the
extra smoothness conditions on its first and second derivatives only in
a neighborhood of $\Gamma_s$.
Moreover, under the assumptions of Theorem~\ref{th:ratiolimit}, it is
possible to calculate an upper bound for the error $E^*(n)$ as a
function of the number of evaluations performed by AEP.
Indeed, \eqref{eq:final4} can be rewritten as
%
\begin{equation}\label{eq:compute1}
E^*(n) \leq A f_*(d)^n.
\end{equation}
We now denote by $M(n)$ the total number of evaluations of the joint
distribution $H$ performed by AEP after the $n$th iteration. Then,
$M(n)$ (as well as the computational time used) is proportional to the
number of simplexes $f_S(d)^{n-1}$ passed to the $n$th iteration. For
all $n \geq2$, we have that
\begin{eqnarray}\label{eq:compute2}
M(n) &=& \sum_{k=0}^{n-1} 2^d f_S(d)^k= \frac{2^d}{f_S(d)-1} \bigl(f_S(d)^n-1\bigr)\nonumber
\\[-8pt]\\[-8pt]
&\geq& \biggl(\frac{2^d}{f_S(d)-1}-1 \biggr) f_S(d)^{n}= B f_S(d)^{n}.\nonumber
\end{eqnarray}
Here, $B$ is a positive constant depending only on the dimension $d$.
Combining~\eqref{eq:compute1} and~\eqref{eq:compute2} gives
%
\begin{equation}\label{eq:cerror}
E^*(n) \leq A \biggl( \frac{M(n)}{B} \biggr)^{\sfrac{\ln f_*(d) }{\ln f_S(d) }}.
\end{equation}
Then,~\eqref{eq:cerror} provides an upper bound on the AEP
approximation error $E^*(n)$ as a function of the number of evaluations
performed. The polynomial rate of\vspace*{1pt} convergence $\frac{\ln f_*(d) }{\ln
f_S(d) }$ of this bound depends only on the dimensionality $d$. In
Table~\ref{ta:simplexes}, we calculate this bound for dimensions $d
\leq8$. These numbers can be useful in order to compare the efficiency
of AEP with that of other algorithms, such as Monte Carlo methods (see
Section~\ref{se:comparison} and Table~\ref{ta:QRM}).



\section{Applications} \label{se:applications}

In this section, we test the AEP algorithm on some risk vectors
$(X_1,\ldots,X_d)$ of financial and actuarial interest.
For illustrative reasons, we will provide the joint distribution
function $H$ in terms of the marginal distributions $F_{X_i}$ and a
copula $C$. For the theory of copulas, we refer the reader to~\citep{rN06}.

In Table~\ref{ta:1A}, we consider a two-dimensional portfolio ($d=2$)
with Pareto marginals, that is,
\[
F_{X_i}(x)=\PP[X_i \leq x]=1-(1+x)^{-\theta_i},\qquad x \geq0, i=1,2,
\]
with tail parameters $\theta_1=0.9$ and $\theta_2=1.8$.
We couple these Pareto marginals via a Clayton copula $C=C^{Cl}_{\delta
}$ with
\begin{eqnarray*}
C^{Cl}_{\delta}(u_1,\ldots,u_d)= ( u_1^{-\delta}+u_2^{-\delta}+\cdots
+u_d^{-\delta}-d+1 )^{-1/\delta},\qquad u_k \in[0,1], k=1,\ldots,d.
\end{eqnarray*}
The parameter $\delta$ is set to 1.2.
For the portfolio described above, we compute the approximation
$P_n(s)$ (see~\eqref{eq:defpnd}) at some given thresholds $s$ and for
different numbers of iterations $n$ of the algorithm.
The thresholds $s$ are chosen in order to have estimates in the center
as well as in the (heavy) tail of the distribution.
For each $n$, we provide the computational time needed to obtain the
estimate on an Apple MacBook (2.4 GHz Intel Core 2 Duo, 2~GB RAM). Of
course, computational times may vary depending on the hardware used for
computations. We also provide the estimates obtained by using the
estimator $P_n^*(s)$, as defined in~\eqref{eq:estformula2}.

\begin{figure}

\includegraphics{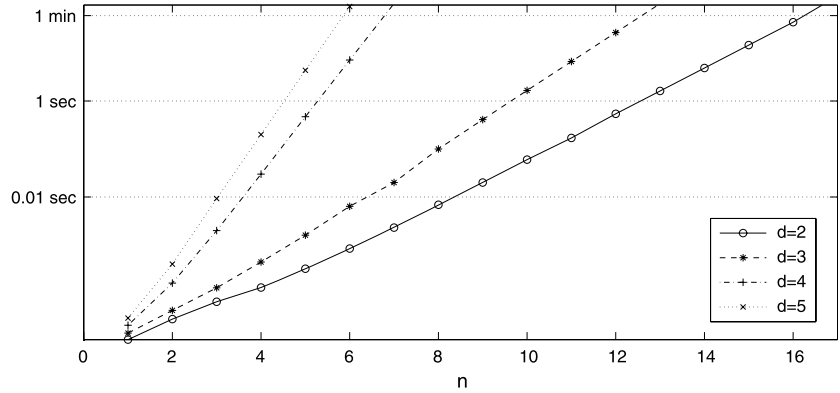}

\caption{AEP computation time (on a log-scale) as a function of the
number of iterations $n$, for dimensions $ 2 \leq d \leq5$.}\label{fi:comptimes}
\end{figure}

For all iterations $n$ and thresholds $s$, in Table~\ref{ta:1A}, we
provide the differences $P_n(s)-P_{16}(s)$ or $P_n^*(s)-P_{16}(s)$.
This has been done in order to show the speed of convergence of the
algorithm and the increase in accuracy due to extrapolation. The choice
of $n=16$ as the reference value in Table~\ref{ta:1A} represents the
maximum number of iterations allowed by the memory (2~GB RAM) of our laptop.
However, for a two-dimensional vector, we see that all iterations after
the seventh leave the first eight decimal digits of the probability
estimate unaltered for all the thresholds. Thus, the estimate $P_7(s)$
(0.01 seconds) could already be considered reasonably accurate. We also
note that, on average, extrapolation allows the accuracy of the
estimates to be increased by two decimal digits without increasing
computational time.

In Tables~\ref{ta:1B} ($d=3$) to~\ref{ta:1D} ($d=5$) we perform the
same analysis for different Clayton--Pareto models in which we
progressively increase the number of random variables used.
In Tables~\ref{ta:1B}--\ref{ta:1D}, the numbers $n=13$ for $d=3$, $n=7$
for $d=4$ and $n=6$ for $d=5$ again represent the maximum number of
iterations allowed by the memory (2~GB RAM) of our laptop.

AEP shows good convergence results for all dimensions $d$ and
thresholds $s$ under study.
In higher dimensions $d$, the extrapolation technique still seems to
provide some relevant extra accuracy.
Memory constraints made estimates for $d \geq6$ prohibitive. For
dimensions $2 \leq d \leq5$, Figure~\ref{fi:comptimes} shows that the
average computational time needed by AEP to provide a single estimate
increases exponentially in the number of iterations $n$.
These average computational times have been computed based on several
portfolios of Pareto marginals coupled by a Clayton copula.

\begin{sidewaystable}
\tablewidth=\textwidth
\caption{Values for $P_n(s)$ and $P_n^*(s)$ (starred columns) for the
sum of two Pareto distributions with parameters $\theta_1=0.9$ and
$\theta_2=1.8$, coupled by a Clayton copula with parameter $\delta
=1.2$; for all $n<16$, we give the difference from the reference value
$P_{16}(s)$}\label{ta:1A}
\begin{tabular*}{\tablewidth}{@{\extracolsep{4in minus 4in}}llllllll@{}}
\hline
& $n = 16$ &$n = 7$ &$n = 7^*$ &$n = 10$& $n = 10^*$&$n = 13$ & $n =13^*$ \\
& (reference value, 49.25 s)& (0.01 s) & (0.01 s) & (0.06 s) & (0.06 s) & (1.61 s) & (1.61 s) \\
\hline
$s=10^0$& $0.315835041363441$ & $-4.46e{-}09 $ & $ -1.46e{-}11 $ & $ -6.16e{-}12 $ & $ -3.70e{-}14 $ & $ -3.97e{-}14 $ & $ -2.95e{-}14 $ \\
$s=10^2$ & $ 0.983690398913354 $ & $ -3.10e{-}10 $ & $ +1.83e{-}09 $ & $ -1.85e{-}12 $ & $ -5.68e{-}13 $ & $ -6.64e{-}13 $ & $ -6.96e{-}13 $ \\
$s=10^4$ & $ 0.999748719229367 $ & $ -6.62e{-}08 $ & $ -4.13e{-}08 $ & $ -6.41e{-}12 $ & $ +6.38e{-}11 $ & $ -1.24e{-}12 $ & $ -1.26e{-}12 $ \\
$s=10^6$ & $ 0.999996018908404 $ & $ -1.63e{-}09 $ & $ -1.22e{-}09 $ & $ -5.40e{-}11 $ & $ -3.89e{-}11 $ & $ -7.80e{-}13 $ & $ -5.07e{-}13 $ \\
\hline
\end{tabular*}\vspace*{12pt}
%
%
\caption{This is the same as Table~\protect\ref{ta:1A}, but for the sum of
three Pareto distributions with parameters $\theta_1=0.9$, $\theta
_2=1.8$ and $\theta_3=2.6$, coupled by a Clayton copula with parameter
$\delta=0.4$}\label{ta:1B}
\begin{tabular*}{\tablewidth}{@{\extracolsep{4in minus 4in}}llllllll@{}}
\hline
& $n = 13$ &$n = 7$ &$n = 7^*$&$n = 9$& $n = 9^*$&$n = 11$& $n = 11^*$
\\
& (reference value, 118.50 s)& (0.02 s) & (0.02 s) & (0.41 s) & (0.41 s) & (6.65 s) & (6.65 s) \\ \hline
$s=10^0$ & $ 0.190859309689430 $ & $ -2.28e{-}06 $ & $ +8.80e{-}07 $ & $ -8.53e{-}08 $ & $ +3.31e{-}08 $ & $ -3.15e{-}09 $ & $ +1.32e{-}09 $ \\
$s=10^2$ & $ 0.983659549676444 $ & $ -1.76e{-}05 $ & $ +1.13e{-}06 $ & $ -6.55e{-}07 $ & $ +3.01e{-}07 $ & $ -2.17e{-}08 $ & $ +1.11e{-}08 $ \\
$s=10^4 $ & $ 0.999748708770280 $ & $ -1.72e{-}06 $ & $ -1.12e{-}06 $ & $ -3.86e{-}07 $ & $ -2.39e{-}07 $ & $ -6.43e{-}08 $ & $ -2.95e{-}08 $ \\
$s=10^6$ & $ 0.999996018515584 $ & $ -2.78e{-}08 $ & $ -1.83e{-}08 $ & $ -6.61e{-}09 $ & $ -4.26e{-}09 $ & $ -1.35e{-}09 $ & $ -7.66e{-}10 $ \\
\hline
\end{tabular*}\vspace*{12pt}
%
%
\caption{This is the same as Table~\protect\ref{ta:1A}, but for the sum of four
Pareto distributions with parameters $\theta_1=0.9$, $\theta_2=1.8$,
$\theta_3=2.6$ and $\theta_4=3.3$, coupled by a Clayton copula with
parameter $\delta=0.2$}\label{ta:1C}
\begin{tabular*}{\tablewidth}{@{\extracolsep{4in minus 4in}}llllllll@{}}
\hline
& $n = 7$ &$n = 4$&$n = 4^*$&$n = 5$& $n = 5^*$ &$n = 6$& $n = 6^*$ \\
& (reference value, 107.70 s)& (0.03 s) & (0.03 s) & (0.47 s) & (0.47 s) & (7.15 s) & (7.15 s) \\ \hline
$s=10^1$ & $ 0.833447516734442 $ & $ -6.31e{-}03 $ & $ +9.42e{-}05 $ & $ -2.21e{-}03 $ & $ +3.71e{-}04 $ & $ -6.04e{-}04 $ & $ +4.00e{-}04 $ \\
$s=10^2$ & $ 0.983412214152579 $ & $ -1.61e{-}03 $ & $ -4.95e{-}04 $ & $ -7.14e{-}04 $ & $ -1.54e{-}04 $ & $ -2.45e{-}04 $ & $ +5.01e{-}05 $\\
$s=10^3 $ & $ 0.997950264030106 $ & $ -2.14e{-}04 $ & $ -7.37e{-}05 $ & $ -9.91e{-}05 $ & $ -2.70e{-}05 $ & $ -3.60e{-}05 $ & $ +3.68e{-}06 $\\
$s=10^4$ & $ 0.999742266243751 $ & $ -2.69e{-}05 $ & $ -9.30e{-}06 $ & $ -1.25e{-}05 $ & $ -3.42e{-}06 $ & $ -4.54e{-}06 $ & $+4.52e{-}07 $\\
\hline
\end{tabular*}
\end{sidewaystable}

\begin{figure}[b]

\includegraphics{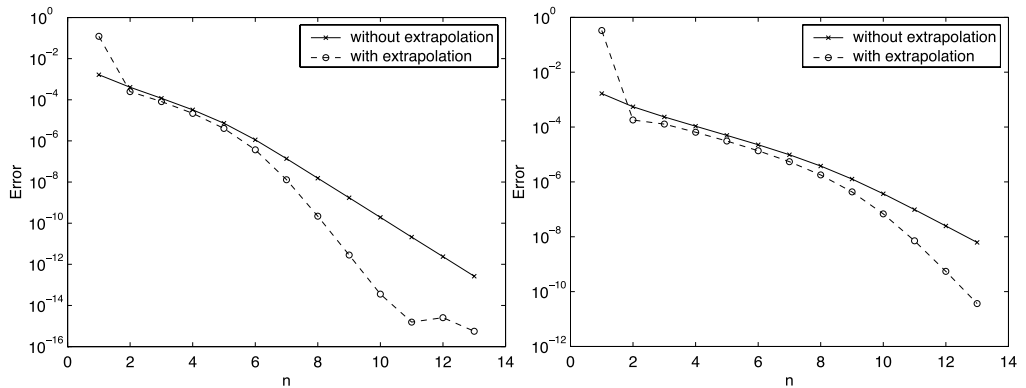}

\caption{Error from the AEP algorithm with and without the use of the
extrapolation technique for two test portfolios: two (left) and three
(right) independent Pareto marginals with parameters $\theta_i=i$, $i=1,2,3.$}
\label{fi:extrapolation}
\end{figure}

Note that Tables~\ref{ta:1A}--\ref{ta:1D} provide information about the
convergence of the algorithm to a certain value, but do not say
anything about the correctness of the limit. Indeed, we do not have
analytical methods to compute
$V_H[\I(\veccc{0},s)]$ when the vector $(X_1,\ldots,X_d)$ has a general
dependence structure (copula) $C$.

In practice, it is possible to test the accuracy of AEP in particular
cases when the $X_i$ are independent or comonotonic.
Some test cases are analyzed in Tables~\ref{ta:2A} ($d=2$) to~\ref{ta:2C} ($d=4$), where
we still assume that we have Pareto marginals, but coupled by a Gumbel
copula $C=C^{Gu}_{\gamma}$, in which the parameter $\gamma\geq1$ is
allowed to vary.
Formally, for $u_k \in(0,1]$, $k=1,\ldots,d$, we have
\begin{eqnarray*}
C^{Gu}_{\gamma}(u_1,\ldots,u_d)=\exp \bigl( - [(-\ln u_1)^\gamma+ (-\ln u_2)^\gamma+\cdots+ (-\ln u_d)^\gamma ]^{1/ \gamma} \bigr).
\end{eqnarray*}
In the tables mentioned above, the multivariate model varies from
independence ($\gamma=1$) to comonotonicity ($\gamma=+\infty$). In
these two extreme (with respect to the dependence parameter~$\gamma$)
cases, we compare the analytical values for $V_H[\I(\veccc{0},s)]$ with
their AEP estimates.
Tables~\ref{ta:1A}--\ref{ta:1D} show that the extrapolated estimator
$P_n^*(s)$ provides accurate estimates within a very reasonable
computational time. A comparison with alternative methods is discussed
in Section~\ref{se:comparison}.

The possibility of computing the value $V_H [\I(\veccc{0},s) ]$
independently from AEP also allows us to test more specifically the
effect of extrapolation.
For this, we consider two- and three-dimensional vectors of independent
Pareto marginals. Figure~\ref{fi:extrapolation} shows the increase of
accuracy due to extrapolation. Therefore, under a smooth model for $H$
(see Theorem~\ref{th:ratiolimit}), the extrapolated estimator
$P_n^*(s)$ is to be preferred over $P_n(s)$.

\begin{sidewaystable}
\tablewidth=\textwidth
\caption{This is the same as Table~\protect\ref{ta:1A}, but for the sum of five
Pareto distributions with parameters $\theta_1=0.9$, $\theta_2=1.8$,
$\theta_3=2.6$, $\theta_4=3.3$ and $\theta_5=4$, coupled by a Clayton
copula with parameter $\delta=0.3$}\label{ta:1D}
\begin{tabular*}{\tablewidth}{@{\extracolsep{4in minus 4in}}llllllll@{}}
\hline
& $n = 6$ &$n = 3$&$n = 3^*$ &$n = 4$& $n = 4^*$ &$n = 5$ & $n = 5^*$\\
& (reference value, 92.91 s)& (0.01 s) & (0.01 s) & (0.20 s) & (0.20 s) & (4.37 s) & (4.37 s) \\ \hline
$s=10^1$ & $ 0.824132635126808 $ & $ -3.12e{-}02 $ & $ +3.89e{-}03 $ & $ -1.55e{-}02 $ & $ +5.66e{-}04 $ & $ -7.77e{-}03 $ & $ +1.46e{-}04 $ \\
$s=10^2$ & $ 0.983253494805448 $ & $ -5.30e{-}03 $ & $ +5.07e{-}05 $ & $ -2.86e{-}03 $ & $ -3.57e{-}04 $ & $ -1.54e{-}03 $ & $ -1.90e{-}04 $ \\
$s=10^3 $ & $ 0.997930730055234 $ & $ -6.72e{-}04 $ & $ -5.23e{-}06 $ & $ -3.66e{-}04 $ & $ -5.29e{-}05 $ & $ -1.99e{-}04 $ & $ -2.83e{-}05 $\\
$s=10^4$ & $ 0.999739803851201 $ & $ -8.45e{-}05 $ & $ -7.22e{-}07 $ & $ -4.61e{-}05 $ & $ -6.67e{-}06 $ & $ -2.51e{-}05 $ & $-3.57e{-}06 $\\
\hline
\end{tabular*}\vspace*{12pt}
\caption{Values for $P_
{n}^*(s)$ for the sum of two Pareto distributions with parameters
$\theta_i=i$, $i=1,2,$ coupled by a Gumbel copula with parameter~$\gamma
$; the values in the first and last columns are calculated
analytically; the computational time for each estimate in this table is
0.53 seconds with $n=12$}\label{ta:2A}
\begin{tabular*}{\tablewidth}{@{\extracolsep{4in minus 4in}}llllllll@{}}
\hline
& $\gamma=1$ (exact) & $\gamma=1$& $\gamma= 1.25$& $\gamma= 1.5$& $\gamma= 1.75$&$\gamma= +\infty$ & $\gamma= +\infty$ (exact) \\
\hline
$s=10^0 $ & $ 0.2862004 $ & $ 0.2862004 $ & $ 0.3280000 $ & $ 0.3527174 $ & $ 0.3682522 $ & $ 0.4108029 $ & $ 0.4108027$ \\
$s=10^2 $ & $ 0.9898913 $ & $ 0.9898913 $ & $ 0.9895957 $ & $ 0.9894472 $ & $ 0.9893640 $ & $ 0.9891761 $ & $ 0.9891761$ \\
$s=10^3 $ & $ 0.9989990 $ & $ 0.9989990 $ & $ 0.9989857 $ & $ 0.9989798 $ & $ 0.9989766 $ & $ 0.9989700 $ & $ 0.9989700 $\\
$s=10^4 $ & $ 0.9999000 $ & $ 0.9999000 $ & $ 0.9998995 $ & $ 0.9998993$ & $ 0.9998992 $ & $ 0.9998990 $ & $ 0.9998990 $\\ \hline
\end{tabular*}\vspace*{12pt}
\caption{This is the same as Table~\protect\ref{ta:2A}, but for the sum of
three Pareto distributions with parameters $\theta_i=i$, $i=1,2,3,$
coupled by a Gumbel copula with parameter $\gamma$; the computational
time for each estimate in this table is 6.65 seconds with $n=11$}\label{ta:2B}
\begin{tabular*}{\tablewidth}{@{\extracolsep{4in minus 4in}}llllllll@{}}
\hline
& $\gamma=1$ (exact) & $\gamma=1$& $\gamma= 1.25$& $\gamma= 1.5$&
$\gamma= 1.75$& $\gamma= +\infty$ & $\gamma= +\infty$ (exact) \\
\hline
$s=10^1 $ & $ 0.1709337 $ & $ 0.1709337 $ & $ 0.2348582 $ & $ 0.2743918 $ & $ 0.2994054 $ & $ 0.3667285 $ & $ 0.3666755$\\
$s= 10^2 $ & $ 0.9898380 $ & $ 0.9898380 $ & $ 0.9893953 $ & $ 0.9891754 $ & $ 0.9890526 $ & $ 0.9887811 $ & $ 0.9887760$\\
$s=10^3 $ & $ 0.9989985 $ & $ 0.9989985 $ & $ 0.9989812 $ & $ 0.9989734 $ & $ 0.9989692 $ & $ 0.9989604 $ & $ 0.9989606$\\
$s= 10^4 $ & $ 0.9999000 $ & $ 0.9999000 $ & $ 0.9998994 $ & $ 0.9998992 $ & $ 0.9998991 $ & $ 0.9998988 $ & $ 0.9998988$ \\
\hline
\end{tabular*}
\end{sidewaystable}

Of course, the AEP algorithm can be used to find estimates for the
quantile function, that is, for the inverse of the distribution of the
sum $S_d$. Such quantiles are especially useful in finance and
insurance, where they are generally referred to as value-at-risk (VaR)
or return periods.
In Table~\ref{ta:4A}, we calculate, by numerical inversion, VaR at
different quantile levels $\alpha$ for two different three-dimensional
portfolios of risks. In order to calculate VaR values, we use
root-finding algorithms like the bisection method.

We finally note that the choices of copula families (Clayton, Gumbel)
and marginal distributions used in this section are purely illustrative
and do not in any way affect the functioning of the AEP algorithm. The
same performances were reached for vectors showing negative dependence,
as in the case of $d$ Pareto marginals coupled by a Frank copula with
negative parameter.

\begin{table}
\tabcolsep=0pt
\caption{This is the same as Table~\protect\ref{ta:2A}, but for the sum of four
Pareto distributions with parameters $\theta_i=i$, $i=1,2,3,4$ coupled
by a Gumbel copula with parameter $\gamma$; the computational time for
each estimate in this table is 7.15 seconds with $n=6$}\label{ta:2C}
\begin{tabular*}{\tablewidth}{@{\extracolsep{4in minus 4in}}llllllll@{}}
\hline
& $\gamma=1$ (exact) & $\gamma=1$& $\gamma= 1.25$& $\gamma= 1.5$& $\gamma= 1.75$& $\gamma= +\infty$ & $\gamma= +\infty$ (exact) \\
\hline
$s= 10^0 $ & $ 0.1040880 $ & $ 0.1040713 $ & $ 0.1762643 $ & $ 0.2244387 $ & $ 0.2555301 $ & $ 0.3387648 $ & $ 0.3390320$\\
$s=10^2 $ & $ 0.9898032 $ & $ 0.9896608 $ & $ 0.9892592 $ & $ 0.9890502 $ & $ 0.9889268 $ & $ 0.9886415 $ & $ 0.9885287$\\
$s=10^3 $ & $ 0.9989981 $ & $ 0.9989732 $ & $ 0.9989652 $ & $ 0.9989616 $ & $ 0.9989595 $ & $ 0.9989743 $ & $ 0.9989558$\\
$s=10^4 $ & $ 0.9999000 $ & $ 0.9998973 $ & $ 0.9998973 $ & $ 0.9998973 $ & $ 0.9998973 $ & $ 0.9998973 $ & $ 0.9998987$ \\ \hline
\end{tabular*}
\end{table}

\begin{table}[b]
\caption{Value-at-risk for: (a) a three-dimensional portfolio with
marginals $F_1=\operatorname{Exp}(0.2)$, $F_2=\mbox{Logn}(\mu=-0.5,\sigma
^2=9/2)$, $F_3=\operatorname{Pareto}(1.2)$ and a Gumbel copula with $\gamma=1.3$;
(b) a three-dimensional portfolio with Pareto marginals with parameters
$\theta_1=0.8$, $\theta_2=1$, $\theta_3=2$ and a Clayton copula with
$\delta=0.4$; the computation of all VaR estimates needs approximately
49 seconds with $n=10$}\label{ta:4A}
\begin{tabular*}{\textwidth}{@{\extracolsep{4in minus 4in}}d{1.6}d{6.2}d{8.2}@{}}
\hline
\multicolumn{1}{@{}l}{$\alpha$} & \multicolumn{1}{l}{$\mbox{VaR}_{\alpha}^{(a)}$} & \multicolumn{1}{l@{}}{$\mbox{VaR}_{\alpha}^{(b)}$}
\\ \hline
 0.9 & 24.76 & 32.87 \\
 0.99 & 137.67 & 445.36 \\
 0.999 & 700.20 & 6864.58 \\
 0.9999 & 3394.78 & 112442.31 \\
 0.99999 & 17962.78 & 1903698.40 \\
 0.999999 & 108190.96 & 32889360.00 \\ \hline
\end{tabular*}
\end{table}

\begin{figure}

\includegraphics{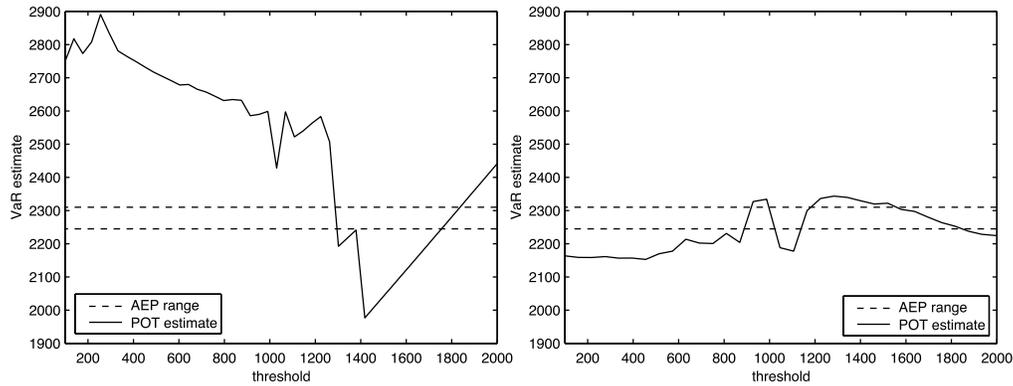}

  \caption{Estimates of $\operatorname{VaR}_{0.999}$ for the sum of the five Pareto marginals
  described in Table~\protect\ref{ta:1D}, as a function of the threshold used for estimation.
  Estimates are obtained via POT from $M=5e{-}03$ (left) and $M=1e{-}06$ (right) simulated data.
  Along with POT estimates, we give the numerical range for the 0.999-quantile obtained via AEP.}\label{fi:pot}
\end{figure}

The accuracy of AEP is not sufficient to estimate high level quantiles
in dimensions $d=4,5,$ as done in Table~\ref{ta:4A} for some
three-dimensional portfolios. The algorithm can, however, be used to
compute a numerical range for the quantiles of the sum of four and five
random variables.
The error resulting from AEP in these higher dimensions turns out to be
extremely small if compared to the error due to statistical inference.
As a comparison to statistical methods, we estimate the VaR of the sum
of the five Pareto marginals described in Table~\ref{ta:1D} via extreme
value theory (EVT) methodology in its ``peaks over threshold'' (POT)
form; see~\citep{MNFE05}, Section 7.2. We set the quantile level $\alpha
=0.999$, a value not uncommon in several risk management applications
in insurance and finance. The POT method is widely used for calculating
quantiles in the presence of heavy-tailed risks and is known to perform
very well in the case of exact Pareto models, such as the one studied
here. In order to focus on the statistical error produced by the POT
method, we use, as data, a sample of $M$ realizations from the
portfolio described in Table~\ref{ta:1D}. It is well known that the
statistical reliability of the POT approach is very sensitive to the
choice of the threshold $u$ beyond which a GPD distribution is fitted.
In Figure~\ref{fi:pot}, we plot the VaR estimates obtained by choosing
different thresholds $u$.
The picture on the left is obtained by generating $M=5000$ data, while
the one on the right uses $M=10^6$ simulations. It is remarkable that,
even in an ideal $10^6$ data world, the statistical range of variation
of the
VaR estimates obtained via POT is broader than the numerical VaR range
calculated via AEP. Moreover, the POT range of values depends on the
specific sample used for estimation, while the AEP range is deterministic.
In the next section, we will compare AEP with more competitive
numerical techniques such as Monte Carlo, quasi-Monte Carlo and
quadrature methods.

\section{A comparison with Monte Carlo, quasi-Monte Carlo and
quadrature methods}\label{se:comparison}

For the estimation of $V_H [{\I(\veccc{0},s)} ]$, the main competitors of
the AEP algorithm are probably Monte Carlo and quasi-Monte Carlo methods.
Given $M$ points $\veccc{x}_1,\ldots,\veccc{x}_M$ in $\I(\veccc{0},s)$, it is
possible to approximate $V_H [{\I(\veccc{0},s)} ]$ by the average of the
density function $v_H$ evaluated at those points, that is,
%
\begin{equation}\label{eq:intd}
V_H [{\I(\veccc{0},s)} ]=\int_{\I(\veccc{0},s)} \ud H(\veccc{x})\simeq\frac
{s^d}{d!} \frac{1}{M} \sum_{i=1}^{M} v_H(\veccc{x}_i).
\end{equation}
If the $\veccc{x}_i$'s are chosen to be (pseudo-)randomly distributed,
this is the \textit{Monte Carlo} (MC) method. If the $\veccc{x}_i$'s are
chosen as elements of a low-discrepancy sequence, this is the \textit{quasi-Monte Carlo} (QMC) method. A \textit{low-discrepancy} sequence
is a totally deterministic sequence of vectors that generates
representative samples from a uniform distribution on a given set. With
respect to Monte Carlo methods, the advantage of using quasi-random
sequences is that points cannot cluster coincidentally on some region
of the set. However, randomization of a low-discrepancy sequence often
improves performance; see~\citep{EL00}.

In recent years, various methods and algorithms have been developed in
order to reduce the variance of MC and QMC estimators and to obtain
probabilities of (rare) events with reasonable precision and effort.
For details on the theory of \textit{rare event simulation} within MC
methods, we refer the reader to~\citep{AG07,pG04,McL05,McL08}.
For an introduction to quasi-Monte Carlo methods and recent
improvements, we refer to, for instance, \citep{hN92} and~\citep{EL00}.
A comprehensive overview of both methods is given in~\citep{SW00}.

Using central limit theorem arguments, it is possible to show that
traditional MC, using (pseudo-)random numbers, has a convergence rate
of $\mathrm{O}(M^{-1/2})$, independently of the number of dimensions $d$.
QMC can be much faster than MC with errors approaching $\mathrm{O}(M^{-1})$ in
optimal cases (see~\citep{wM98}), but the worst theoretic rate of
convergence decreases with the dimension $d$ as $\mathrm{O}((\log M)^d M^{-1})$;
see~\citep{hN92}.
In applications to finance and insurance, it is more common to get
results closer to the best rate of convergence if the density $v_H$ is
smooth, that is, has a Lipschitz-continuous second derivative. In this
case, it is possible to show that the convergence rate is at least
$\mathrm{O}((\log M)^d M^{-3/2})$; see~\citep{CMO97}.
In Table~\ref{ta:QRM}, we compare convergence rates of MC and QMC
methods with respect to the AEP rates (depending on $d$), as provided
in Section~\ref{se:numerics}. We thus expect a well-designed QMC
algorithm to perform better, asymptotically, than AEP under a smooth
probability model and for dimensions $d \geq4$. Because of the
computational issues for AEP in higher dimensions, we restrict our
attention to $d \leq5$ in Table~\ref{ta:QRM}.

\begin{table}
\caption{Asymptotic convergence rates of the AEP, standard MC and QMC methods}\label{ta:QRM}
\begin{tabular*}{\textwidth}{@{\extracolsep{4in minus 4in}}lllll@{}}
\hline
 $d$ & 2 & 3 & 4 & 5  \\
\hline
 AEP (upper bound) & $M^{-3}$ & $M^{-1.5}$ & $M^{-0.54}$ &$M^{-0.34}$  \\
 MC & $M^{-0.5}$ & $M^{-0.5}$ & $M^{-0.5}$ & $M^{-0.5}$  \\
 QMC (best) & $M^{-1}$ & $M^{-1}$ & $M^{-1}$ & $M^{-1}$  \\
 QMC (worst) & $M^{-1}(\log M)^2 $ & $M^{-1}(\log M)^3 $ &$M^{-1}(\log M)^4 $ &$M^{-1}(\log M)^5 $  \\ \hline
\end{tabular*}
\end{table}

Don McLeish kindly adapted an algorithm using a randomized Korobov
low-discrepancy sequence to the portfolio leading to Table~\ref{ta:1A}.
The parameters for the sequence are those recommended in~\citep{HL07}.
The standard errors (s.e.'s) are obtained by independently randomizing
ten (part~(a) of the table) and fifty (part~(b) of the table) sequences
with 1 million terms each, corresponding to $M=1e{-}07$ (a) and $M=5e{-}07$ (b).
The average CPU times are, of course, on a different machine (IBM
Thinkpad 2.5 GHz Intel Core 2 Dual, 4 GB RAM).
In Table~\ref{ta:QRM1}, we provide the comparison between QMC and AEP
extrapolated estimates. The results seem to be coherent with Table~\ref{ta:QRM} above. For the same precision, AEP is much faster than QMC in
the two-dimensional example and slightly slower for $d=4$. Recall that,
in higher dimensions, programming a randomized Korobov rule is much
more demanding than using AEP.

\begin{table}[b]
\begin{subtable}
\caption{AEP and QMC (using Korobov sequence) estimates for $V_H [{\I(\veccc{0},s)} ]$ for the sum of
Two Pareto distributions with
parameters $\theta_1=0.9$ and $\theta_2=1.8$, coupled by a Clayton
copula with parameter $\delta=1.2$}\label{ta:QRM1}
\begin{tabular*}{\tablewidth}{@{\extracolsep{4in minus 4in}}llll@{}}
\hline
$s$ & AEP estimate ($n=14$, 4.87 s) & QMC estimate ($M=$1$e{-}$07, 6.6 s)& QMC s.e. \\
 \hline
 $10^0$ & 0.315835041363413 & 0.3158345 & $+2.7e{-}06$ \\
 $10^2$ & 0.983690398912470 & 0.98369106 & $+1.0e{-}06$ \\
 $10^4$ & 0.999748719228038 & 0.99974872 & $+1.5e{-}07$ \\
 $10^6$ & 0.999996018907752 & 0.999996 & $+4.0e{-}08$ \\
\hline
\end{tabular*}
\caption{AEP and QMC (using Korobov sequence) estimates for $V_H [{\I(\veccc{0},s)} ]$
for the sum of four Pareto distributions with
parameters $\theta_1=0.9$, $\theta_2=1.8$, $\theta_3=2.6$, $\theta
_4=3.3$, coupled by a Clayton copula with parameter $\delta=0.2$;
computational times are also provided}\label{ta:QRM1b}
\begin{tabular*}{\tablewidth}{@{\extracolsep{4in minus 4in}}llll@{}}
\hline
$s$ & AEP estimate ($n=7$, 107.70 s) & QMC estimate ($M=$5$e{-}$07, 95 s) &QMC s.e. \\ \hline
 $10^1$& 0.833826902853978 & 0.83380176 & $+3.6e{-}06$ \\
 $10^2$ & 0.983565803484355 & 0.98362452 & $+9.0e{-}07$ \\
 $10^3$ & 0.997972831330699 & 0.997997715 & $+2.3e{-}07$ \\
 $10^4$ & 0.999745113409911 & 0.999748680 & $+5.0e{-}08$ \\
\hline
\end{tabular*}
\end{subtable}
\end{table}

What is important to stress here is that in MC and randomized QMC
methods similar to the one applied in Table~\ref{ta:QRM1}, the final
estimates contain a source of randomness.
Contrary to this, the AEP algorithm is deterministic, being solely
based on geometric properties of a certain domain.
Moreover, the accuracy of MC and QMC methods is generally lost for
problems in which the density $v_H$ is not smooth or cannot be given in
closed form, and comes at the price of an adaptation of the sampling
algorithm to the specific example under study. Recall that the AEP
algorithm does not require the density of the distribution $H$ in
analytic form, nor does it have to assume overall smoothness. Finally,
the precision of MC methods depends on the threshold $s$ at which $V_H
[{\I(\veccc{0},s)} ]$ is evaluated: estimates in the (far) tail of the
distribution will be less accurate.

The re-tailoring, from example to example, of the rule to be iterated
is also common to other numerical techniques for
the estimation of $V_H [{\I(\veccc{0},s)} ]$ such as \textit{quadrature
methods}; see~\citep{DR84} and~\citep{PTSV07} for a review.
However, in the computation of multi-dimensional integrals, as in~\eqref
{eq:intd}, numerical quadrature rules are typically less efficient than
MC and QMC.

When the random variables $X_1,\ldots,X_d$ are exchangeable and
heavy-tailed, some asymptotic approximations of $V_H [{\I(\veccc{0},s)}
]$ for large $s$ can be found in~\cite{BFG06,LGH05} and
references therein.
It is important to remark that the behavior of AEP is not affected by
the threshold $s$ at which $V_H [{\I(\veccc{0},s)} ]$ is computed, nor by
the tail properties of the marginal distributions $F_{X_i}$. This is
particularly interesting as, under heavy-tailedness, the relative error
of MC and QMC methods increases in the tail of the distribution
function of $S_d$.

We are, of course, aware that a well-designed quadrature rule or a
specific quasi-random sequence might perform better than AEP in a
specific example, with respect to both accuracy and computational effort.
However, AEP provides very accurate estimates of the distribution of
sums up to five dimensions in a reasonable time \textit{without} the
need to adapt to the probabilistic model under study.
AEP can handle, in a uniform way, any joint distribution $H$, possibly
in the form of its copula and marginal distributions.
Because of its ease of use and the very weak assumptions upon which it
is based, AEP offers a competitive tool for the computation of the
distribution function of a sum of up to five random variables. A
Web-based, user-friendly version has been programmed and will
eventually be made available.

\section{Final remarks}\label{se:conclusions}

In this paper, we have introduced the AEP algorithm in order to compute
numerically the distribution function of the sum of $d$ random
variables $X_1,\ldots,X_d$ with given joint distribution~$H$.
The algorithm is mainly based on two assumptions: the random variables
$X_i$ are bounded from below and the distribution $H$ has a bounded
density in a neighborhood of the curve $\Gamma_s$ defined in~\eqref
{eq:gammasdef}. Under this last assumption, the sum $S_d$ has to be
continuous at the threshold $s$ where the distribution is calculated,
that is, $\PP[S_d=s]=0$.
When, instead, $V_H [ \Gamma_s ] >0$, the algorithm may fail to
converge. As an example, take two random variables $X_1$ and $X_2$ with
$\PP[X_1=1/2]=\PP[X_2=1/2]=1$. Then, $V_H [ \I(\veccc{0},1) ]=1$, but the
sequence $P_n(1)$ alternates between $0$ and $1$. Similar examples for
arbitrary dimension $d$ can easily be constructed.

If $H$ has at least a bounded density near $\Gamma_s$, then the
convergence of the sequence $P_n(s)$ to the value $V_H [\I(\veccc{0},s)
]$ is guaranteed. As already remarked, the speed of convergence may
vary, depending on the probability mass of a neighborhood of~$\Gamma
_s$. Tools to increase the efficiency of the algorithm are therefore
much needed in these latter cases.

The AEP algorithm has been shown to converge when $d\leq5$ if the
joint distribution $H$ of the vector $(X_1,\ldots,X_d)$ has a bounded
density $v_H$. Under some extra smoothness assumptions on $v_H$,
convergence holds when $d\leq8$.
All of these conditions can be weakened to hold only in a neighborhood
of the curve $\Gamma_s$ and are satisfied by most examples which are
relevant in practice.

We were not able to prove convergence of AEP in arbitrary dimensions,
although we conjecture this to hold. The main problem in higher
dimensions is the non-monotonicity of $P_n(s)$ and $P_n^*(s)$. This
results from the fact that the $s_n^k$'s, as defined in~\eqref
{eq:defofs}, may be positive as well as negative.
From a geometric point of view, the main problem is the fact that the
simplexes $\I_{n+1}^{k}$, $k=1,\ldots,N^n$, passed to the $(n+1)$th
iteration of the algorithm, are generally not disjoint for $d>2$. As
illustrated in Table~\ref{ta:volfact}, the sum of the Lebesgue measures
of the $\I_{n+1}^{k}$'s is increasing in the number $n$ of iterations
when $d>6$, while their union always lies in some neighborhood of the
curve $\Gamma_s$.
A general convergence theorem may need a volume decomposition different
from~\eqref{eq:decd} and using only a family of \textit{disjoint}
simplexes, or else an extension of the extrapolation technique.

We also remark that the statement of a general convergence theorem will
not entail any practical improvement of AEP, since memory constraints
limit the use of the algorithm to dimension $d \leq5$.
However, in these manageable dimensions, we expect the AEP convergence
rates to be better than their upper bounds given in Table~\ref{ta:simplexes}.

Apart from the study of convergence of AEP in higher dimensions, in
future research, we will also address an extension of the algorithm to
more general aggregating functions $\psi(X_1,\ldots,X_d)$ and the study
of an adaptive (i.e., depending on $H$) and more efficient (in terms of
new simplexes produced at each iteration) decomposition of the simplexes.

\begin{appendix}
\section*{\texorpdfstring{Appendix: Proof of~\protect\eqref{eq:decd}}{Appendix: Proof of (3.1) }}
\label{se:appendixa}

Recall that, in Section~\ref{se:intro}, we denoted
by $\veccc{i}_0,\ldots, \veccc{i}_{N}$ all of the $2^d$ vectors in $\{0,1\}
^d$, with $\veccc{i}_0=(0,\ldots,0)$, $\veccc{i}_k=\veccc{e}_k$, $k=1,\ldots
,d$, and $\veccc{i}_{N}=\veccc{1}=(1,\ldots,1)$, where $N=2^{d}-1$. Also,
recall that $\# \veccc{i}$ denotes the number of $1$'s in the vector $\veccc
{i}$, for instance, $\# \veccc{i}_0=0$, $\# \veccc{i}_N=d$.
\begin{theorem}\label{ap:mainth}
For any $\veccc{b}\in\mathbb{R}^d$, $h\in\mathbb{R}$ and $\alpha\in
[1/d,1)$, we have that
\begin{eqnarray*}
V_H [\I(\veccc{b},h) ] = V_H [\Q(\veccc{b},\alpha h) ]
+\sum_{j=1}^N m^j V_H [\I(\veccc{b}^j,h^j) ],
\end{eqnarray*}
where, for all $j=1,\ldots,N$,
\begin{eqnarray}\label{eq:appdef}
\veccc{b}^{j}&=&\veccc{b}+\alpha h \veccc{i}_j,\qquad
h^{j}=(1- \# \veccc{i}_j \alpha)h,\nonumber
\\[-8pt]\\[-8pt]
m^j &=&\cases{
(-1)^{1+\#\veccc{i}_j},         &\quad \mbox{if } $\#\veccc{i}_j < 1/\alpha$,\cr
0,                              &\quad \mbox{if } $\#\veccc{i}_j = 1/\alpha$,\cr
(-1)^{d+1-\#\veccc{i}_j},       &\quad \mbox{if } $\#\veccc{i}_j > 1/\alpha$.
}\nonumber
\end{eqnarray}
\end{theorem}

Note that~\eqref{eq:appdef} is equivalent to~\eqref{eq:decd} under the
notation introduced in Section~\ref{se:algo}.
In order to prove the above theorem, we need some lemmas. In the
following, $\delta_{ij}$ denotes the Kronecker delta, that is,
\[
\delta_{ij}=\cases{
0, &\quad\mbox{if }$i \neq j,$\cr
1, &\quad\mbox{if }$i =j$.
}
\]
\begin{lemma}\label{le:ap1}
Fix $i,j \in D$ with $i \neq j$. Then, for any $h,s \in\mathbb{R}$
with $hs \geq0$ and $\veccc{b}\in\mathbb{R}^d$, we have that
\[
\I(\veccc{b}+h\veccc{e}_i,s) \cap\I(\veccc{b}+h\veccc{e}_j,s)
= \cases{
\I(\veccc{b}+h\veccc{e}_j+h\veccc{e}_i,s-h), &\quad \mbox{if } $| h | < | s |$,
\cr
\varnothing, &\quad \mbox{if } $| h | \geq | s |$.
}
\]
\end{lemma}
\begin{pf}\vspace*{-13pt}
\begin{pf*}{Proof of $ \subset$}
First, assume $0<s\leq h$. By definition~\eqref{de:simplex}, for a
vector $\veccc{x}\in\I(\veccc{b}+h\veccc{e}_i,s)$, we have that
\begin{eqnarray*}
x_k>b_k+\delta_{ik}h, k \in D\quad \mbox{and}\quad \sum_{k=1}^d (x_k-b_k- \delta
_{ik}h) \leq s,
\end{eqnarray*}
from which it follows that
\[
x_j \leq b_j+s - \sum_{k \neq j} (x_k-b_k-\delta_{ik}h) < b_j +s \leq b_j+h,
\]
that is, $\veccc{x} \notin\I(\veccc{b}+h\veccc{e}_j,s)$.
Now, assume that $0<h < s$. For a vector $\veccc{x} \in\I(\veccc{b}+h\veccc
{e}_i,s) \cap\I(\veccc{b}+h\veccc{e}_j,s)$, we have that
%
\begin{equation}\label{eq:lapp1}
x_k-b_k>0, k \in D \qquad\mbox{with } x_i>b_i+h \mbox{ and } x_j>b_j+h.
\end{equation}
Again, $\veccc{x} \in\I(\veccc{b}+h\veccc{e}_i,s)$, therefore $\sum_{k=1}^d
(x_k-(b_k+h\delta_{ik})) \leq s$.
Subtracting $h$ from both sides of the last inequality, we obtain
%
\begin{equation}\label{eq:lapp2}
\sum_{k=1}^d \bigl(x_k-(b_k+h\delta_{ik}+h\delta_{jk}) \bigr) \leq s-h.
\end{equation}
Equations~\eqref{eq:lapp1} and~\eqref{eq:lapp2} show that $\veccc{x} \in
\I(\veccc{b}+h\veccc{e}_j+h\veccc{e}_i,s-h)$. The case $h,s<0$ is
analogous.\noqed
\end{pf*}

\begin{pf*}{Proof of $\supset$}
If $0<s \leq h$, there is nothing to show. Suppose, then, that $0<h<s$.
For any fixed $\veccc{x}\in\I(\veccc{b} + h\veccc{e}_j + h\veccc{e}_i,s-h)$,
\eqref{eq:lapp2} holds with $x_k-(b_k+h\delta_{ik}+h\delta_{jk})>0$, $k
\in D$. By adding $h \delta_{jk}$ in the sum on the left-hand side and
$h$ to the right-hand side of~\eqref{eq:lapp2}, we find that
%
\begin{equation}\label{eq:lapp2bis}
\sum_{k=1}^d \bigl(x_k-(b_k+h\delta_{ik}) \bigr)\leq s.
\end{equation}
Since $ (x_k-(b_k+h\delta_{ik}) )$ is still positive for all $k \in D$,
\eqref{eq:lapp2bis} shows that $\veccc{x} \in\I(\veccc{b}+h\veccc{e}_i,s)$.
By similar reasoning, we also have that $\veccc{x} \in\I(\veccc{b}+h\veccc{e}_j,s)$.
The case $h,s<0$ is analogous; the case $hs=0$ is trivial.
\end{pf*}\noqed
\end{pf}

\begin{lemma}\label{le:ap2}
For any $\veccc{b}\in\mathbb{R}^d$, $h\in\mathbb{R}$ and $\alpha\in
(0,1)$, we have that
\[
\I(\veccc{b},h)\setminus\Q(\veccc{b},\alpha h)=
\bigcup_{k=1}^d \I(\veccc{b}+\alpha h \veccc{e}_k,h-\alpha h).
\]
\end{lemma}

\begin{pf}\vspace*{-13pt}
\begin{pf*}{Proof of $ \subset$}
First, assume that $h>0$. If $\veccc{x} \in\I(\veccc{b},h) \setminus\Q(\veccc
{b},\alpha h)$, then $x_k > b_k$, $k \in D$ and $\sum_{k=1}^d (x_k -
b_k) \leq h$, while, by definition~\eqref{de:cube}, there exists a $j
\in D$ such that $x_j -b_j> \alpha h$. For this $j$, it is then
possible to write
%
\begin{eqnarray}\label{eq:lapp3}
\sum_{k=1}^d \bigl(x_k - (b_k+\delta_{jk}\alpha h)\bigr) \leq h-\alpha h
\qquad\mbox{with } x_k - (b_k+\delta_{jk}\alpha h)>0, k \in D,
\end{eqnarray}
which yields $\veccc{x}\in\I(\veccc{b}+\alpha h \veccc{e}_j, h - \alpha h)
\subset\bigcup_{k=1}^d\I(\veccc{b}+\alpha h \veccc{e}_k,h-\alpha h).$
\noqed
\end{pf*}

\begin{pf*}{Proof of $ \supset$}
Let $\veccc{x} \in\bigcup_{k=1}^d \I(\veccc
{b}+\alpha h \veccc{e}_k,h-\alpha h)$, meaning that there exists $j \in
D$ for which $\veccc{x}$ satisfies~\eqref{eq:lapp3}. It follows that $x_j
> b_j+\alpha h$ (hence $\veccc{x} \notin\Q(\veccc{b},\alpha h)$) and $\sum
_{k=1}^d (x_k - b_k) \leq h-\alpha h+\alpha h=h$. Noting that~\eqref
{eq:lapp3} also implies $x_k>b_k$, $k \in D$, we finally obtain that
$\veccc{x} \in\I(\veccc{b},h)\setminus\Q(\veccc{b},\alpha h)$.
The case $h<0$ is analogous, while the case $h=0$ is trivial.
\end{pf*}\noqed
\end{pf}
\begin{lemma}\label{le:ap3}
For any $\veccc{b}\in\mathbb{R}^d$, $h\in\mathbb{R}$ and $\alpha\in
[1/d,1)$, we have that
\[
\Q(\veccc{b},\alpha h) \setminus\I(\veccc{b},h)= \I(\veccc{b}+\alpha h \veccc
{1},h-\alpha d h)\cap\Q(\veccc{b},\alpha h).
\]
\end{lemma}
\begin{pf}\vspace*{-13pt}
\begin{pf*}{Proof of $ \subset$}
If $\alpha=1/d$, then the lemma is straightforward. So, choose $\alpha
\in(1/d,1)$ and assume $h>0$. If $\veccc{x} \in\Q(\veccc{b},\alpha h)
\setminus\I(\veccc{b},h)$, then $x_k>b_k$ for all $k \in D$. Since $\veccc
{x} \notin\I(\veccc{b},h)$, it follows that $\sum_{i=1}^d (x_i-b_i)>h$.
Since $x_k \leq b_k+\alpha h$ for all $k \in D$, we can write
%
\begin{equation}\label{eq:lapp4}
\sum_{k=1}^d (x_k-b_k-\alpha h)>h-\alpha dh \qquad\mbox{with }
x_k-(b_k+\alpha h) \leq0 \mbox{ for all } k \in D.
\end{equation}
As $h-d\alpha h=h(1-d\alpha)<0$, we conclude that
$\veccc{x} \in\I(\veccc{b}+\alpha h \veccc{1},h-\alpha d h)$ and, hence, by
assumption,
$\veccc{x} \in\I(\veccc{b}+\alpha h \veccc{1},h-\alpha d h)\cap\Q(\veccc
{b},\alpha h)$.
\end{pf*}

\begin{pf*}{Proof of $ \supset$} Let $\veccc{x} \in \I(\veccc{b}+\alpha h \veccc
{1},h-\alpha d h)\cap\Q(\veccc{b},\alpha h)$. Due to $h-\alpha dh<0$, it
follows that~\eqref{eq:lapp4} holds, implying that $\sum_{k=1}^d
(x_k-b_k)>h$, that is, $\veccc{x} \notin\I(\veccc{b},h)$. The case $h<0$
is analogous, while the case $h=0$ is trivial.
\end{pf*}\noqed
\end{pf}

We are now ready to prove the main result in this appendix.
\begin{pf*}{Proof of Theorem~\protect\ref{ap:mainth}}
The case $h=0$ is
trivial. Suppose, then, that $h\neq0$. From the general property of
two sets $A,B$ that $B= (A \cup(B \setminus A)) \setminus(A\setminus
B)$, $(A\setminus B) \subset A \cup(B \setminus A)$ and $A \cap(B
\setminus A)=\varnothing$, it follows that
%
\begin{eqnarray} \label{eq:lapp0}
V_H [\I(\veccc{b},h) ] = V_H [\Q(\veccc{b},\alpha h) ]
+ V_H [\I(\veccc{b},h) \setminus\Q(\veccc{b},\alpha h) ]
- V_H [\Q(\veccc{b}, \alpha h) \setminus\I(\veccc{b}, h) ].
\end{eqnarray}
Using the notation $\I^k=\I(\veccc{b}+\alpha h \veccc{e}_k,h-\alpha h)$,
Lemma~\ref{le:ap2} implies, for the second summand in~\eqref{eq:lapp0}, that
%
\begin{eqnarray}\label{eq:lapp6}
V_H [\I(\veccc{b},h) \setminus\Q(\veccc{b},\alpha h) ] =
V_H \Biggl[ \bigcup_{k=1}^d \I^k \Biggr]=\sum_{k=1}^d (-1)^{k+1}
\sum_{I\subset D, | I | = k} V_H \biggl[ \bigcap_{i\in I} \I^i \biggr].
\end{eqnarray}
Fixing $I\subset D$ with $I=\{n_1,\ldots,n_k\}$,
iteratively using Lemma~\ref{le:ap1} yields
\begin{eqnarray*}
\bigcap_{i\in I} \I(\veccc{b}+\alpha h \veccc{e}_{n_i},h-\alpha h) =
\cases{
\I\Biggl(\veccc{b}+\alpha h\displaystyle  \sum_{j=1}^k \veccc{e}_{n_j},h(1- k \alpha)\Biggr), &\quad\mbox{if }$k \alpha<1$,\cr
\varnothing,                                                        &\quad\mbox{if }$k \alpha\geq1$.
}
\end{eqnarray*}
Substituting this last expression into~\eqref{eq:lapp6} implies
that
\begin{eqnarray}\label{eq:lapp7}
V_H [\I(\veccc{b},h) \setminus\Q(\veccc{b},\alpha h) ]
&=&\mathop{\sum_{k\in D,}}_{k \alpha<1} (-1)^{k+1}
\mathop{\sum_{\veccc{i}_r \in\{0,1\}^d, }}_{ \#\veccc{i}_r=k}
V_H \bigl[\I\bigl(\veccc{b}+\alpha h \veccc{i}_r,h(1- k \alpha)\bigr) \bigr] \nonumber\\[-8pt]\\[-8pt]
&=&\mathop{\sum_{\veccc{i} \in\{0,1\}^d, }}_{ 0<\#\veccc{i}<1/ \alpha}
(-1)^{\# \veccc{i}+1} V_H \bigl[\I\bigl(\veccc{b}+\alpha h \veccc{i},h(1- \# \veccc{i}
\alpha)\bigr) \bigr].\nonumber
\end{eqnarray}
Using Lemma~\ref{le:ap3} for the third summand in~\eqref{eq:lapp0}, we
can also write that\vspace*{-1pt}
%
\begin{eqnarray}\label{eq:lapp9}
&&V_H [\Q(\veccc{b}, \alpha h) \setminus\I(\veccc{b}, h) ] \nonumber\\
&&\quad= V_H
[\I(\veccc{b}+\alpha h \veccc{1},h-\alpha d h)\cap\Q(\veccc{b},\alpha h) ]
\\
&&\quad= V_H [\I(\veccc{b}+\alpha h \veccc{1},h-\alpha d h) ] - V_H [\I(\veccc
{b}+\alpha h \veccc{1},h-\alpha d h)\setminus\Q(\veccc{b},\alpha h) ].\nonumber
\end{eqnarray}
Note that if $\alpha=1/d$, then the quantity in~\eqref{eq:lapp9}
is zero. We can hence assume that $\alpha\neq1/d$.
Observing that $\Q(\veccc{b},\alpha h) = \Q(\veccc{b}+\alpha h \veccc
{1},-\alpha h)$ and defining $\hat{\veccc{b}}=\veccc{b}+\alpha h \veccc{1}$,
$\hat{\alpha}=-\alpha/(1-\alpha d)>1/d$ and $\hat{h}=h(1-\alpha d)$, we
can write\vspace*{-1pt}
\[
V_H [\I(\veccc{b}+\alpha h \veccc{1},h-\alpha d h)\setminus\Q(\veccc
{b},\alpha h) ]=V_H [\I(\hat{\veccc{b}},\hat{h})\setminus\Q(\hat{\veccc
{b}},\hat{\alpha}\hat{h}) ].
\]
Note that the right-hand side of the previous equation is empty if $\hat
{\alpha} \geq1$, that is, $\alpha\in(1/d, 1/(d-1)]$.
At this point, equation~\eqref{eq:lapp7} yields\vspace*{-1pt}
\begin{eqnarray*}
&&V_H [\I(\veccc{b}+\alpha h \veccc{1},h-\alpha d h)\setminus\Q(\veccc
{b},\alpha h) ]\\
&&\quad=
\mathop{\sum_{\veccc{i} \in\{0,1\}^d, }}_{ 0<\#\veccc{i}<1/ \hat{\alpha}} (-1)^{\# \veccc{i}+1} V_H \bigl[\I\bigl(\hat{\veccc{b}}+\hat
{\alpha} \hat{h} \veccc{i}, \hat{h}(1- \# \veccc{i} \hat{\alpha})\bigr) \bigr]\\
&&\quad=\mathop{\sum_{\veccc{i} \in\{0,1\}^d,}}_{0<\#\veccc{i}<d-1/\alpha}
(-1)^{\# \veccc{i}+1} V_H \bigl[\I\bigl(\veccc{b}+\alpha h (\veccc{1}-\veccc{i}),h\bigl(1-\alpha(d -\# \veccc{i})\bigr)\bigr) \bigr].
\end{eqnarray*}
Substituting $\hat{\veccc{i}}=\veccc{1}-\veccc{i}$ ($\#\hat{\veccc{i}}=d- \#
\veccc{i}$) into the previous equation, we can equivalently write\looseness=1\vspace*{-1pt}
\begin{eqnarray}\label{eq:lapp8}
&&V_H [\I(\veccc{b}+\alpha h \veccc{1},h-\alpha d h)\setminus\Q(\veccc
{b},\alpha h) ]\nonumber
\\[-8pt]\\[-8pt]
&&\quad=
\mathop{\sum_{\hat{\veccc{i}} \in\{0,1\}^d, }}_{ 1/\alpha<\#\hat{\veccc{i}}<d} (-1)^{d- \# \hat{\veccc{i}}+1} V_H \bigl[\I\bigl(\veccc
{b}+\alpha h \hat{\veccc{i}}, h(1-\# \hat{\veccc{i}} \alpha)\bigr) \bigr].\nonumber
\end{eqnarray}
In keeping with what was noted above, this last equation is null in the
aforementioned case in which $\hat{\alpha} \geq1$.
Recalling~\eqref{eq:lapp9} and noting that
\begin{eqnarray*}
\I(\veccc{b}+\alpha h \veccc{1},h-\alpha d h)
= \I\bigl(\veccc{b}+\alpha h \veccc{i}_N,h(1-\# \veccc{i}_N \alpha)\bigr),
\end{eqnarray*}
we obtain
\vspace{-7pt}
\begin{eqnarray}\label{eq:lapp10}
&&V_H [\Q(\veccc{b},h) \setminus\I(\veccc{b},\alpha h) ]\nonumber\\
&&\quad=
V_H \bigl[\I\bigl(\veccc{b}+\alpha h \veccc{i}_N,h(1- \# \veccc{i}_N \alpha)\bigr) \bigr]\nonumber
\\[-8pt]\\[-8pt]
&&\qquad{}-
\mathop{\sum_{\hat{\veccc{i}} \in\{0,1\}^d, }}_{ 1/ \alpha<\#\hat{\veccc{i}}<d} (-1)^{d- \# \hat{\veccc{i}}+1}
V_H \bigl[\I\bigl(\veccc{b}+\alpha h \hat{\veccc{i}}, h(1-\# \hat{\veccc{i}} \alpha)\bigr)
\bigr]\nonumber\\[-4pt]
&&\quad=\mathop{\sum_{\hat{\veccc{i}} \in\{0,1\}^d, }}_{ 1/ \alpha<\#\hat{\veccc
{i}} \leq d} (-1)^{d- \# \hat{\veccc{i}}} V_H \bigl[\I\bigl(\veccc{b}+\alpha h \hat
{\veccc{i}}, h(1-\# \hat{\veccc{i}} \alpha)\bigr) \bigr].\nonumber
\end{eqnarray}
Finally, recalling the definitions in~\eqref{eq:appdef}, we substitute
equations~\eqref{eq:lapp7} and~\eqref{eq:lapp10} into~\eqref{eq:lapp0}
to obtain
\vspace{-4pt}
\begin{eqnarray*}
V_H [\I(\veccc{b},h) ] &=& V_H [\Q(\veccc{b},\alpha h) ]
+ \mathop{\sum_{\veccc{i} \in\{0,1\}^d, }}_{ 0<\#\veccc{i}<1/ \alpha}
(-1)^{\# \veccc{i}+1} V_H \bigl[\I\bigl(\veccc{b}+\alpha h \veccc{i},h(1- \# \veccc{i}\alpha)\bigr) \bigr]
\\
&&{}- \mathop{\sum_{\hat{\veccc{i}} \in\{0,1\}^d, }}_{ 1/ \alpha<\#\hat{\veccc
{i}} \leq d} (-1)^{d- \# \hat{\veccc{i}}} V_H \bigl[\I\bigl(\veccc{b}+\alpha h \hat
{\veccc{i}}, h(1-\# \hat{\veccc{i}} \alpha)\bigr) \bigr] \\
&=&V_H [\Q(\veccc{b},\alpha h) ]
+\sum_{j=1}^{N} m^j V_H [\I(\veccc{b}^j,h^j) ].
\end{eqnarray*}
\upqed
\end{pf*}
\end{appendix}

\section*{Acknowledgements}
The authors are grateful to Don McLeish for providing relevant
comments on the paper and the example illustrated in Table~\ref{ta:QRM}.
Giovanni Puccetti would like to thank RiskLab and the
Forschungsinstitut f\"ur Mathematik (FIM) of the Department of
Mathematics, ETH
Z\"urich, for its financial support and kind hospitality. Philipp
Arbenz would like to thank SCOR for financial support toward the final
stages of writing this paper. The final version of the paper was
written while Paul Embrechts was visiting the Institute for
Mathematical Sciences at the National University of Singapore.
Finally, the authors would like to thank two anonymous referees and an
Associate Editor for several valuable comments which significantly
improved the paper.

\printhistory

\end{document}